\documentclass[11pt]{article}

\usepackage{amsmath,latexsym,amssymb,amsfonts,amsbsy, amsthm}
\usepackage{graphicx}
\usepackage{empheq}
\usepackage{color}
\usepackage{ulem}
\def\inte#1{
	\displaystyle\mathop{#1\kern0pt}^\circ }

\numberwithin{equation}{section}

\numberwithin{Proposition}{section}
\newtheorem{Theorem}{Theorem}
\numberwithin{Theorem}{section}
\newtheorem{Lemma}{Lemma}
\numberwithin{Lemma}{section}

\numberwithin{Remark}{section}

\numberwithin{Corollary}{section}
\date{}

\begin{document}
\title{Stability of Couette flow for 2D Boussinesq system in a uniform magnetic field}

\author{Dongfen Bian \footnote{ School of Mathematics and Statistics,
		Beijing Institute of Technology, Beijing 100081, China. Email: {\tt biandongfen@bit.edu.cn/dongfen\_bian@brown.edu}.} \and
Shouyi Dai\footnote{School of Mathematics and Statistics,
	Beijing Institute of Technology, Beijing 100081, China. Email: {\tt daishoui@outlook.com}.} \and
	Jingjing Mao \footnote{School of Mathematics and Statistics,
		Beijing Institute of Technology, Beijing 100081, China. Email: {\tt mao.jingjing@outlook.com}.} }

\maketitle
\begin{abstract}
In this paper, we consider the Boussinesq equations with magnetohydrodynamics convection in the domain $\mathbb{T} \times \mathbb{R}$ and establishes the nonlinear stability of
	 the Couette flow $(\mathbf{u}_{sh} = (y,0), \mathbf{b}_{sh} = (1,0), p_{sh} = 0, \theta_{sh} = 0$). The novelty in this paper is that we design a new Fourier multiplier operator by using the properties of the enhanced dissipation to overcome the difficult term $\partial_{xy}(-\Delta)^{-1}j$ in the linearized and nonlinear system. Then, we prove the asymptotic  stability for the linearized system.
Finally, we establish the nonlinear stability for the full system by bootstrap principle. 
\end{abstract}

\begin{center}
	\begin{minipage}{120mm}
	\end{minipage}
\end{center}
\begin{center}
	\begin{minipage}{120mm}
	\end{minipage}
\end{center}

\section{Introduction and Main Results}
In this paper, we consider the following 2D incompressible Boussinesq equations for magnetohydrodynamics (MHD) convection \cite{bian2020nonlinear,majda2003introduction}
\begin{equation}\label{system_full_dissipation}
\begin{cases}
\mathbf{u}_{t} + (\mathbf{u}\cdot \nabla)\mathbf{u} - (\mathbf{b} \cdot \nabla) \mathbf{b} + \nabla p - \nu \Delta \mathbf{u} = \theta \mathbf{e}_{2} , \\
\mathbf{b}_{t} + (\mathbf{u}\cdot \nabla)\mathbf{b} - (\mathbf{b} \cdot \nabla) \mathbf{u} - \mu \Delta \mathbf{b} = 0 , \\
\theta_{t} + (\mathbf{u}\cdot \nabla)\theta - \eta\Delta\theta = 0 , \\
\nabla \cdot \mathbf{u} = \nabla \cdot \mathbf{b} = 0 .
\end{cases}
\end{equation}
The unknowns are the velocity field  $\mathbf{u} = (u^{1},u^{2})$, the magnetic field  $\mathbf{b} = (b^{1},b^{2})$, the temperature $\theta$ and the scalar pressure $p$. In addition, we denote here by $\nu$ the fluid viscosity,  $\mu$ the magnetic diffusivity, $\eta$ the thermal diffusivity and  $e_{2}=(0, 1)$. The spatial domain $\Omega$ here is taken to be $\Omega = \mathbb{T} \times \mathbb{R}$ with $\mathbb{T} = [0,2\pi]$ being the periodic box and $\mathbb{R}$ being the whole line.

Physically, the first equation of  (\ref{system_full_dissipation}) represents the
conservation law of the momentum with the effect of the buoyancy $\theta e_{2}$. The second equation of  (\ref{system_full_dissipation}) shows that the electromagnetic field is governed by the Maxwell equation. The third  equation is a balance of the temperature convection and diffusion. For more physics and numerical simulations, the interested readers may refer to \cite{bian2020nonlinear,kulikovskiy1965magnetohydrodynamics} and the references therein. 
Two fundamental problems, the global regularity problem and the stability problem, have been among the main driving forces in advancing the mathematical theory on the  Boussinesq-MHD system. 
Significant progress has been made on the global regularity of the  nonlinear Boussinesq-MHD system
 \cite{bian20162,bian2017initial,bian2020global,bian2016initial,larios2017local,li2019global,liu2019global, yu2020global,zhai2018global}. The goal of this paper is  the nonlinear stability around the Couette flow  $(\mathbf{u}_{sh} = (y,0), \mathbf{b}_{sh} = (1,0), p_{sh} = 0, \theta_{sh} = 0$).

 When the fluid is not affected by the temperature, that is, $\theta \equiv 0$, then the equations (\ref{system_full_dissipation}) become the incompressible MHD system and govern the dynamics of the velocity and the magnetic field in electrically conducting fluids such as plasmas and reflect the basic physics conservation laws. Also the issue of stability on the MHD equations has been extensively studied. Liss considered the sobolev stability threshold of 3D Couette flow in a uniform magnetic field \cite{liss2020sobolev}. In \cite{ji2020resistive}, the author  obtained  the stability and large-time behavior of perturbations near a stationary solution of the 2D resistive MHD equation.  Further background and motivation for the MHD system may be found in \cite{desjardins1998remarks,
sermange1983some} and references therein.

When  the fluid is not affected by the Lorentz force, that is, $\mathbf{b} \equiv 0$, then the equations (\ref{system_full_dissipation}) become the classical Boussinesq system. The Boussinesq system reflects the basic physics laws obeyed by buoyancy-driven fluids. It is one of the most frequently used model for atmospheric and oceanographic flows and serves as the centerpiece in the study of the Rayleigh convection \cite{constantin1996heat, Doering1995Applied, majda2003introduction, pedlosky2013geophysical}.    
Important progress has been made on the stability and large-time behavior \cite{castro2019asymptotic,doering2018long,masmoudi2020stability, tao2020stability,yang2018linear,Zill20,Zill2011}. Very recently, Deng, Wu and  Zhang prove the nonlinear stability of Couette flow for the 2D Boussinesq system \cite{deng2020stability}.  Motivated by this work,  we  consider the Boussinesq system in the presence of magnetic field.

Note that when  the fluid is not affected by  the temperature and the Lorentz force, that is, $\theta\equiv 0$ and $\mathbf{b} \equiv 0$, then the equations (\ref{system_full_dissipation}) become the Navier-Stokes equations, and become the Euler equations if without viscosity. For the Couette flow of Euler and Navier-Stokes equations, there are many interesting results 
\cite{bedrossian2015dynamics,
bedrossian1506dynamics,
bedrossian2017stability,
bedrossian2019stability,
BM13,
BMV16,
chen-li-wei-zhang,
chen-wei-zhang,
Reddy,
Rom73,
transition_threshold_wei_d_and_zhang_z}.

In this paper we mainly focus on the nonlinear stability  for the full system (\ref{system_full_dissipation})  and consider the  Couette flow 
\[
  \mathbf{u}_{sh} = (y,0),\qquad \mathbf{b}_{sh} = (1,0), \qquad p_{sh} = 0, \qquad \theta_{sh} = 0.
\]
The perturbations around this Couette flow take the form of 
\[
  \mathbf{\tilde{u}} = (\tilde{u}^{1}, \tilde{u}^{2}) = (u^{1} - y,u^{2}) ,\qquad \mathbf{\tilde{b}} = (\tilde{b}^{1}, \tilde{b}^{2}) = (b^{1} - 1,b^{2}),\qquad \tilde{p} = p,\qquad \tilde{\theta} = \theta,
\]
and satisfy the following system
\begin{equation*}\label{system_perturbation_with_full_dissipation}
\begin{cases}
\tilde{u}^{1}_{t} + (\mathbf{\tilde{u}}\cdot \nabla)\tilde{u}^{1} + y\partial_{x}\tilde{u}^{1} + \tilde{u}^{2} - (\mathbf{\tilde{b}}\cdot \nabla)\tilde{b}^{1} - \partial_{x} \tilde{b}^{1} + \partial_{x} \tilde{p} - \nu \Delta \tilde{u}^{1} = 0, \\
\tilde{u}^{2}_{t} + (\mathbf{\tilde{u}}\cdot \nabla)\tilde{u}^{2} + y\partial_{x}\tilde{u}^{2}  - (\mathbf{\tilde{b}}\cdot \nabla)\tilde{b}^{2} - \partial_{x} \tilde{b}^{2} + \partial_{y} \tilde{p} - \nu \Delta \tilde{u}^{2} = \tilde{\theta}, \\
\tilde{b}^{1}_{t} + (\mathbf{\tilde{u}}\cdot \nabla)\tilde{b}^{1} + y\partial_{x}\tilde{b}^{1}  - (\mathbf{\tilde{b}}\cdot \nabla)\tilde{u}^{1} - \partial_{x} \tilde{u}^{1} - \tilde{b}^{2} - \mu \Delta \tilde{b}^{1} = 0, \\
\tilde{b}^{2}_{t} + (\mathbf{\tilde{u}}\cdot \nabla)\tilde{b}^{2} + y\partial_{x}\tilde{b}^{2}  - (\mathbf{\tilde{b}}\cdot \nabla)\tilde{u}^{2} - \partial_{x} \tilde{u}^{2} - \mu \Delta \tilde{b}^{2} = 0, \\
\tilde{\theta}_{t} + (\mathbf{\tilde{u}}\cdot \nabla)\tilde{\theta} + y\partial_{x} \tilde{\theta} - \eta\Delta\tilde{\theta} = 0, \\
\nabla \cdot \mathbf{\tilde{u}} = \nabla \cdot \mathbf{\tilde{b}} = 0.
\end{cases}
\end{equation*}
The corresponding perturbed vorticity and current density near the steady vorticity  $w_{sh} = -1$ and the steady current density $j_{sh} = 0$ take the form of 
\[
  \tilde{w} = \partial_{x} \tilde{u}^{2} - \partial_{y} \tilde{u}^{1},\qquad \tilde{j} = \partial_{x} \tilde{b}^{2} - \partial_{y} \tilde{b}^{1}, 
\]
and verify, together with $\tilde{\theta}$, 
the following system 
\begin{equation}\label{nonlinear_stability_system}
\begin{cases}
  \partial_{t}\tilde{w} + (\mathbf{\tilde{u}} \cdot \nabla) \tilde{w} - (\mathbf{\tilde{b} \cdot \nabla})\tilde{j} + y\partial_{x}\tilde{w} - \partial_{x} \tilde{j} - \nu \Delta \tilde{w} = \partial_{x} \tilde{\theta}, \\
  \partial_{t}\tilde{j} + (\mathbf{\tilde{u}} \cdot \nabla) \tilde{j}  - (\mathbf{\tilde{b}} \cdot \nabla) \tilde{w} + y\partial_{x}\tilde{j} - \partial_{x} \tilde{w} - \mu \Delta \tilde{j}  - 2\partial_{x}\tilde{b}^{1} - Q(\nabla\mathbf{\tilde{u}}, \nabla \mathbf{\tilde{b}}) = 0, \\
  \partial_{t}\tilde{\theta} + (\mathbf{\tilde{u}}\cdot \nabla)\tilde{\theta} + y\partial_{x} \tilde{\theta} - \eta\Delta\tilde{\theta} = 0, \\
  \tilde{\mathbf{u}} = -\nabla^{\bot}(-\Delta)^{-1}\tilde{w}, \\
  \tilde{\mathbf{b}} = -\nabla^{\bot}(-\Delta)^{-1}\tilde{j},
\end{cases}
\end{equation}
where $Q(\nabla\mathbf{\tilde{u}}, \nabla \mathbf{\tilde{b}}) = 2\partial_{x}\tilde{b}^{1}(\partial_{x}\tilde{u}^{2} + \partial_{y}\tilde{u}^{1}) - 2\partial_{x}\tilde{u}^{1}(\partial_{x}\tilde{b}^{2} + \partial_{y}\tilde{b}^{1})$.

For notational convenience, we shall write $w$ for $\tilde{w}$, $\theta$ for $\tilde{\theta}$, $j$ for $\tilde{j}$, $\mathbf{u}$ for $\mathbf{\tilde{u}}$ and $\mathbf{b}$ for $\mathbf{\tilde{b}}$ from now on. We rewrite the nonlinear  system (\ref{nonlinear_stability_system}) as follows
\begin{equation}\label{nonlinear_stability_system_rewrite}
\begin{cases}
  \partial_{t}w + (\mathbf{u} \cdot \nabla) w - (\mathbf{b \cdot \nabla})j + y\partial_{x}w - \partial_{x} j - \nu \Delta w = \partial_{x} \theta, \\
  \partial_{t}j + (\mathbf{u} \cdot \nabla) j  - (\mathbf{b} \cdot \nabla) w + y\partial_{x}j - \partial_{x} w - \mu \Delta j  - 2\partial_{x}b^{1} - Q(\nabla\mathbf{u}, \nabla \mathbf{b}) = 0, \\
  \partial_{t}\theta + (\mathbf{u}\cdot \nabla)\theta + y\partial_{x} \theta - \eta\Delta\theta = 0, \\
  \mathbf{u} = -\nabla^{\bot}(-\Delta)^{-1}w, \\
  \mathbf{b} = -\nabla^{\bot}(-\Delta)^{-1}j,
\end{cases}
\end{equation}
where $Q(\nabla\mathbf{u}, \nabla \mathbf{b}) = 2\partial_{x}b^{1}(\partial_{x}u^{2} + \partial_{y}u^{1}) - 2\partial_{x}u^{1}(\partial_{x}b^{2} + \partial_{y}b^{1})$.

The linearization of (\ref{nonlinear_stability_system_rewrite}) takes the form of
\begin{equation}\label{linear_system}
\begin{cases}
  \partial_{t}w + y\partial_{x}w - \partial_{x} j  -\nu \Delta w - \partial_{x} \theta=0, \\
  \partial_{t}j + y\partial_{x}j - \partial_{x} w - \mu \Delta j
    - 2\partial_{x}b^{1} =0, \\
  \partial_{t}\theta + y\partial_{x} \theta - \eta\Delta \theta=0,  \\
  \mathbf{u} = -\nabla^{\bot}(-\Delta)^{-1}w, \\
  \mathbf{b} = -\nabla^{\bot}(-\Delta)^{-1}j,	\\
  w|_{t=0} = w(0),\qquad j|_{t=0} = j(0), \qquad \theta|_{t=0} = \theta(0).
\end{cases}
\end{equation}

Now we state our  main results. 
To make the statement precise, we define, for $f=f(x,y) $ with $(x,y)\in \mathbb{T}\times\mathbb{R}$ and $k\in\mathbb{Z}$,
\begin{equation*}
f_{k}(y):=\frac{1}{2\pi}\int_{\mathbb{T}}f(x,y)e^{-ikx}dx.
\end{equation*}
Then, we write
$f_{\neq}(x, y)=f(x,y)-f_{0}(y)$. 
In addition, for two functions $f$ and $g$ and a norm $\|\cdot \|_{X}$ we write
\begin{equation*}
\|(f,g) \|_{X}=\sqrt{\|f \|_{X}^{2}+\|g\|_{X}^{2}}.
\end{equation*}

For the lineared system \eqref{linear_system}, we can get the following linear stability.
\begin{Theorem}\label{the_linear_stability_result}
  Let $0 < \nu = \mu \leq \eta \leq 1$ and $(w,j,\theta)$ be the solution to the system (\ref{linear_system}) with initial data $(w(0),j(0),\theta(0))$. Then there exist two positive constants $c>0$ and $ C>0$ such that for any $k \in \mathbb{Z}$, $t>0$,
  \begin{equation*}
       \Vert \theta_{k}(t) \Vert_{L^{2}_{y}} \leq C \Vert \theta_{k}(0) \Vert_{L^{2}_{y}} e^{-c\eta^{\frac{1}{3}}\vert k \vert^{\frac{2}{3}}t}, 
      \end{equation*} 
    \begin{equation} \label{the_linear_stability_result1}  
    \Vert\big( w_{k}(t), j_{k}(t)\big) \Vert_{L_{y}^{2}} \leq C \left( \nu^{-2} \Vert\big( w_{k}(0),j_{k}(0)\big) \Vert_{L_{y}^{2}} + \nu^{-6}(\nu^{-1} \vert k \vert)^{\frac{1}{3}}\Vert \theta_{k}(0) \Vert_{L_{y}^{2}} \right)e^{-c \nu^{\frac{13}{3}}\vert k \vert^{\frac{2}{3}}t}.
  \end{equation}
  Moreover, for $N > 0$, there exist $c_{N} > 0$ and $C_{N} > 0$ such that for any $k \in \mathbb{Z}$, $ t>0$,
   \begin{equation}\label{the_linear_stability_result2}
 \Vert D_{y}^{N}\theta_{k}(t) \Vert_{L_{y}^{2}}  \leq  C_{N} e^{-c_{N}\eta^{\frac{1}{3}}\vert k \vert^{\frac{2}{3}}t}\left( \Vert D_{y}^{N}\theta_{k}(0) \Vert_{L^{2}_{y}} + (\eta^{-1}\vert k \vert)^{\frac{N}{3}} \Vert \theta_{k}(0) \Vert_{L^{2}_{y}} \right), 
\end{equation}
\begin{equation}\label{the_linear_stability_result3}
\begin{split}
	&\Vert\big  (D^{N}_{y}w_{k}(t),  D^{N}_{y}j_{k}(t)\big) \Vert_{L_{y}^{2}}  
\\
&\leq  C_{N}
e^{-c_{N}\nu^{\frac{13}{3}}\vert k \vert^{\frac{2}{3}}t}
\Big( \nu^{-2}\Vert\left( D^{N}_{y}w_{k}(0),D^{N}_{y}j_{k}(0) \right)\Vert_{L_{y}^{2}}  
+ \nu^{-6}(\nu^{-1}\vert k \vert)^{\frac{1}{3}}\Vert D^{N}_{y}\theta_{k}(0) \Vert_{L_{y}^{2}}  
 \\&\qquad+ \nu^{-6N}(\nu^{-1}\vert k \vert)^{\frac{N}{3}}
  \big( \nu^{-2}\Vert\big( w_{k}(0), j_{k}(0)\big) \Vert_{L_{y}^{2}} 
 + \nu^{-6}(\nu^{-1}\vert k \vert)^{\frac{1}{3}}\Vert \theta_{k}(0) \Vert_{L_{y}^{2}}\big) \Big) .
    \end{split}
  \end{equation}
\end{Theorem}

The linear stability result in Theorem \ref{the_linear_stability_result} can be converted into  the estimate in physical space by introducing the time-dependent operator, for $t \geq 0$, 
\[
\Lambda_{t}^{b} = \left( 1 - \partial_{x}^{2} - (\partial_{y} + t \partial_{x})^{2} \right)^{\frac{b}{2}},
\]
or, in terms of its symbol, $\Lambda_{t}^{b}(k,\xi) = (1+k^{2}+(\xi + tk)^{2})^{\frac{b}{2}}$. And for any $b\in \mathbb{R}$, it is easy to check that $\Lambda_{t}^{b}$ commutes with $\partial_{t} + y\partial_{x}$. The estimate in physical space is stated in the following theorem.
\begin{Theorem}\label{the_physics_linear_stability_result}
  Let  $(w,j,\theta)$ be the solution to (\ref{linear_system}) with initial data $(w(0),j(0),\theta(0))$, and $0 < \nu = \mu \leq \eta \leq 1$. There exists a constant $C>0$ such that for $b\in \mathbb{R}$,
  \begin{align*}
 &    \Vert \big(\Lambda_{t}^{b}w, \Lambda_{t}^{b}j \big)\Vert_{L_{t}^{\infty}(L^{2})} 
    + \nu^{\frac{1}{2}}\Vert \big(\nabla \Lambda_{t}^{b}w,\nabla\Lambda_{t}^{b}j\big) \Vert_{L_{t}^{2}(L^{2})}
     + \nu^{\frac{1}{6}}\Vert\big( \vert D_{x}\vert^{\frac{1}{3}} \Lambda_{t}^{b}w , \vert D_{x}\vert^{\frac{1}{3}} \Lambda_{t}^{b}j \big)\Vert_{L_{t}^{2}(L^{2})} \\
    &\qquad + \nu^{-4}(\nu\eta)^{-\frac{1}{6}}\left( \Vert \vert D_{x} \vert^{\frac{1}{3}}\Lambda_{t}^{b}\theta \Vert_{L_{t}^{\infty}(L^{2})} + \eta^{\frac{1}{2}}\Vert \nabla \vert D_{x} \vert^{\frac{1}{3}}\Lambda_{t}^{b}\theta \Vert_{L_{t}^{2}(L^{2})} + \eta^{\frac{1}{6}}\Vert\vert D_{x} \vert^{\frac{2}{3}} \Lambda_{t}^{b}\theta \Vert_{L_{t}^{2}(L^{2})} \right) \\
    & \leq C\left( \nu^{-2}\Vert\big( w(0), j(0)\big) \Vert_{H^{b}} + \nu^{-4}(\nu\eta)^{-\frac{1}{6}} \Vert \vert D_{x} \vert^{\frac{1}{3}}\theta(0) \Vert_{H^{b}} \right).
  \end{align*}
\end{Theorem}

The nonlinear stability of the  full  system \eqref{nonlinear_stability_system_rewrite} is stated as follows.
\begin{Theorem}\label{nolinear_theorem}
Assume $\nu = \mu = \eta$, $0 < \nu \leq 1$, $b > 1$, $\beta \geq \frac{11}{2}$, $\delta \geq \beta + \frac{13}{3}$, $\alpha \geq \delta - \beta + \frac{14}{3}$ and assume that the initial data $(w(0),j(0),\theta(0))$ satisfies
\[
  \Vert \theta(0) \Vert_{H^{b}} \leq \varepsilon \nu^{\alpha},
  \qquad \Vert \big(w(0),j(0)\big) \Vert_{H^{b}} \leq \varepsilon \nu^{\beta},
  \qquad \Vert \vert D_{x} \vert^{\frac{1}{3}} \theta(0) \Vert_{H^{b}} \leq \varepsilon \nu^{\delta}, 
\]
for some sufficiently small $\varepsilon > 0$. Then the solution $(w,j,\theta)$ to the system (\ref{nonlinear_stability_system_rewrite}) satisfies 
\begin{align*}
 \Vert \Lambda_{t}^{b}\theta \Vert_{L_{t}^{\infty}(L^{2})} + \nu^{\frac{1}{2}}\Vert \nabla \Lambda_{t}^{b}\theta \Vert_{L_{t}^{2}(L^{2})} + \nu^{\frac{1}{6}}\Vert \vert D_{x} \vert^{\frac{1}{3}} \Lambda_{t}^{b}\theta \Vert_{L_{t}^{2}(L^{2})} 
 + \Vert (-\Delta)^{-\frac{1}{2}}\Lambda_{t}^{b}\theta_{\neq} \Vert_{L_{t}^{2}(L^{2})} \leq C\varepsilon\nu^{\alpha}, 
 \end{align*}
 \begin{align*}
&\Vert \big(\Lambda_{t}^{b}w,\Lambda_{t}^{b}j\big) \Vert_{L_{t}^{\infty}(L^{2})}
   + \nu^{\frac{1}{2}}\Vert \big(\nabla \Lambda_{t}^{b}w,\nabla \Lambda_{t}^{b}j\big) \Vert_{L_{t}^{2}(L^{2})} 
  + \nu^{\frac{1}{6}}\Vert \big(\vert D_{x} \vert^{\frac{1}{3}} \Lambda_{t}^{b}w , \vert D_{x} \vert^{\frac{1}{3}} \Lambda_{t}^{b}j\big)\Vert_{L_{t}^{2}(L^{2})} \\
 &\quad + \Vert\big( (-\Delta)^{-\frac{1}{2}}\Lambda_{t}^{b}w_{\neq},(-\Delta)^{-\frac{1}{2}}\Lambda_{t}^{b}j_{\neq}\big) \Vert_{L_{t}^{2}(L^{2})} \leq C\varepsilon\nu^{\beta},
  \end{align*}
  and
  \begin{align*}
&\Vert \vert D_{x} \vert^{\frac{1}{3}} \Lambda_{t}^{b}\theta \Vert_{L_{t}^{\infty}(L^{2})} + \nu^{\frac{1}{2}}\Vert \nabla \vert D_{x} \vert^{\frac{1}{3}} \Lambda_{t}^{b}\theta \Vert_{L_{t}^{2}(L^{2})} + \nu^{\frac{1}{6}}\Vert \vert D_{x} \vert^{\frac{2}{3}} \Lambda_{t}^{b}\theta \Vert_{L_{t}^{2}(L^{2})} \\
    &\quad+ \Vert (-\Delta)^{-\frac{1}{2}}\vert D_{x} \vert^{\frac{1}{3}}\Lambda_{t}^{b}\theta_{\neq} \Vert_{L_{t}^{2}(L^{2})} \leq C\varepsilon\nu^{\delta},
\end{align*}
where $C$ is a positive constant.
\end{Theorem}

The stability problem of the full system  (\ref{nonlinear_stability_system_rewrite}) is not trivial and more difficult than that in \cite{deng2020stability}. Including the similar difficulty in \cite{deng2020stability}, there is an extra difficult term  $\partial_{xy}(-\Delta)^{-1}j$. More precisely, due to the presence of the buoyancy forcing term, the Sobolev norms or even the $L^{2}$-norm of the velocity field could grow in time if the three linear terms $y\partial_{x}w $, $y\partial_{x}j $ and $y\partial_{x}\theta $ were not included in (\ref{nonlinear_stability_system_rewrite}). Fortunately,  the enhanced dissipation makes the stability of the Couette flow for the system  (\ref{nonlinear_stability_system_rewrite})  possible. 
Following the idea of \cite{deng2020stability}, some operators are designed to extract the properties of enhanced dissipation. 
For the  stablity, we define a function $\phi_{k}$ for $k\neq 0$ as follows,
\begin{equation*}
\phi_{k}(\xi) = 
\begin{cases}
  \frac{6(k^{2} + \xi_{0}^{2})^{2}}{k^{4}} - (2 + \pi) ,\qquad & \xi > 0, \\
  \frac{6(k^{2} + \xi_{0}^{2})^{2}}{(k^{2} + \xi^{2})^{2}} - (2 + \pi), \qquad & \xi \in [-\xi_{0},0], \\
  (4 - \pi) e^{\frac{24\xi_{0}(\xi+\xi_{0})}{(4-\pi)(k^{2}+\xi_{0}^{2})}}, \qquad & \xi \in (-\infty, -\xi_{0}).
\end{cases}
\end{equation*}
Here $\xi_{0}$ is a real positive solution of the equation $\nu \xi_{0} (k^{2} + \xi_{0}^{2}) = 96 \vert k \vert$. It is easy to check that $ \phi_{k} \in C^{1}(\mathbb{R})$.
Then we add $\phi_{k}({\rm sgn}(k)\xi)$ for $k \neq 0$ to the symbol function of the operator constructed in \cite{deng2020stability} to overcome the term $\partial_{xy}(-\Delta)^{-1}j$ and establish the new operator $\mathcal{M}$ for which symbol function satisfies the following properties, 
\begin{align*}
& 2\nu (\xi^{2} + k^{2}) \mathcal{M}(k,\xi) + k\partial_{\xi}\mathcal{M}(k,\xi) \geq \nu(\xi^{2} + k^{2}) + \frac{1}{4}\nu^{\frac{1}{3}}\vert k \vert^{\frac{2}{3}} + \frac{1}{\xi^{2} + k^{2}},  \\
 & 2\nu (\xi^{2} + k^{2}) \mathcal{M}(k,\xi) + k\partial_{\xi}\mathcal{M}(k,\xi) + \mathcal{M}(k,\xi) \frac{4k\xi}{k^{2}+\xi^{2}} \geq \nu(\xi^{2} + k^{2}) + \frac{1}{4}\nu^{\frac{1}{3}}\vert k \vert^{\frac{2}{3}} + \frac{1}{\xi^{2} + k^{2}}.
\end{align*}
The upper bound of the new operator $\mathcal{M}$  is depending on diffusivity. 
This is why we need bigger $\alpha, \beta, \delta$  compared with that in \cite{deng2020stability} in Theorem \ref{nolinear_theorem}.

\section{The linear stability}
This section is devoted to the proofs of the linear stability stated in Theorem \ref{the_linear_stability_result} and Theorem  \ref{the_physics_linear_stability_result}.  In order to prove the desired linear stability results, we construct a new Fourier multiplier operator to overcome the difficult term $\partial_{xy}(-\Delta)^{-1}j$ and the details are in the following subsections. 
\subsection{Proof of Theorem \ref{the_linear_stability_result}}
In this subsection, we  prove Theorem \ref{the_linear_stability_result}.
\begin{proof}
By projecting the equations in (\ref{linear_system}) onto each frequency, we obtain the system in the $y$-variable only,
\begin{equation}\label{linear_system__projecting_y}
\begin{cases}
  \partial_{t}w_{k} + ikyw_{k} - ikj_{k} +   \nu (D^{2}_{y}+ k^{2}) w_{k} = ik \theta_{k}, \\ 
  \partial_{t}j_{k} + ikyj_{k} - ikw_{k}  - 2ikb^{1}_{k} + \mu (D^{2}_{y}+ k^{2}) j_{k} = 0, \\
  \partial_{t}\theta_{k} + iky\theta_{k} +  \eta (D^{2}_{y}+ k^{2}) \theta_{k} = 0, \\
  w_{k}|_{t=0} = w_{k}(0),\quad j_{k}|_{t=0} = \theta_{k}(0), \quad \theta_{k}|_{t=0} = \theta_{k}(0).
\end{cases}
\end{equation}
By taking the $L_{y}^{2}$-inner product with $\theta_{k}$ in the third equation of (\ref{linear_system__projecting_y}), we get
\begin{align*}
 & \langle \partial_{t}\theta_{k},\theta_{k} \rangle_{L_{y}^{2}} + \langle \theta_{k},\partial_{t}\theta_{k} \rangle_{L_{y}^{2}} =  \frac{d}{dt} \langle \theta_{k},\theta_{k} \rangle_{L_{y}^{2}} = \frac{d}{dt} \Vert \theta_{k} \Vert_{L_{y}^{2}}^{2}, \\
& \langle iky\theta_{k},\theta_{k} \rangle_{L_{y}^{2}} + \langle \theta_{k},iky\theta_{k} \rangle_{L_{y}^{2}} = \langle iky\theta_{k},\theta_{k} \rangle_{L_{y}^{2}} + \langle -iky\theta_{k},\theta_{k} \rangle_{L_{y}^{2}} = 0, \\
 & \langle \eta (D^{2}_{y}+ k^{2}) \theta_{k},\theta_{k} \rangle_{L_{y}^{2}} + \langle \theta_{k},\eta (D^{2}_{y}+ k^{2}) \theta_{k} \rangle_{L_{y}^{2}} =  2\eta\Vert D_{y} \theta_{k} \Vert^{2}_{L_{y}^{2}} 
+ 2\eta k^{2}\Vert \theta_{k} \Vert^{2}_{L_{y}^{2}},
\end{align*}
which implies that
\begin{equation}\label{linear_equation_identity1}
 \frac{1}{2}\frac{d}{dt} \Vert \theta_{k} \Vert_{L_{y}^{2}}^{2} + \eta \Vert D_{y} \theta_{k} \Vert^{2}_{L_{y}^{2}} + \eta k^{2}\Vert \theta_{k} \Vert^{2}_{L_{y}^{2}} = 0.
\end{equation}
To further the estimates, we apply  the Fourier multiplier operator defined in \cite{deng2020stability}. If $k\in\mathbb{Z}$ and $k \neq 0$, the multiplier $M_{k}$ is given by
\[
  M_{k}\theta_{k} := \varphi(\eta^{\frac{1}{3}}\vert k \vert^{-\frac{1}{3}}{\rm sgn}(k)D_{y})\theta_{k},
\]
where $\varphi$ is a real-valued, non-decreasing function, and $\varphi \in C^{\infty}(\mathbb{R})$ satisfies  $0 \leq \varphi(x) \leq 1$, $0 \leq \varphi^{'}(x) \leq \frac{1}{4}$ for all $x \in \mathbb{R}$,  and $\varphi^{'}(x) = \frac{1}{4}$ for $ x \in [-1,1]$.  Clearly, $M_{k}$ is a self-adjoint and non-negative Fourier multiplier operator. We take the $L_{y}^{2}$-inner product of the third equation in (\ref{linear_system__projecting_y}) with $M_{k}\theta_{k}$. The following basic identities hold:
\begin{align*}
  2 {\rm  Re} \langle \partial_{t}\theta_{k}, M_{k}\theta_{k} \rangle_{L^{2}_{y}} = &\frac{d}{dt} \langle M_{k}\theta_{k}, \theta_{k} \rangle_{L^{2}_{y}} ,\\
  2 {\rm  Re} \langle \eta(D_{y}^{2} +  k^{2})\theta_{k}, M_{k}\theta_{k} \rangle_{L^{2}_{y}} = & \langle \eta(D_{y}^{2} +  k^{2})\theta_{k}, M_{k}\theta_{k} \rangle_{L^{2}_{y}} + \langle M_{k}\theta_{k}, \eta(D_{y}^{2} +  k^{2})\theta_{k} \rangle_{L^{2}_{y}}\\
  = & \langle \eta(D_{y}^{2} +  k^{2})M_{k}\theta_{k}, \theta_{k} \rangle_{L^{2}_{y}} + \langle \eta(D_{y}^{2} +  k^{2}) M_{k}\theta_{k}, \theta_{k} \rangle_{L^{2}_{y}} \\
  = & \langle 2\eta(D_{y}^{2} +  k^{2})M_{k}\theta_{k}, \theta_{k} \rangle_{L^{2}_{y}} ,\\
  2 {\rm  Re} \langle iky\theta_{k}, M_{k}\theta_{k} \rangle_{L^{2}_{y}} = & \langle [M_{k},iky]\theta_{k},\theta_{k} \rangle_{L^{2}_{y}},
\end{align*} 
where in the last equation we have used the fact that $M_{k}$ is self-adjoint and $iky$ is skew-adjoint. Here the bracket in $[M_{k},iky]$ denotes the standard commutator. Noticing that
\begin{align*}
 \langle [M_{k}, iky]f,g \rangle_{L^{2}_{y}} = & \langle M_{k}f,iky g \rangle_{L^{2}_{y}} + \langle iky f, M_{k}g \rangle_{L^{2}_{y}} \\
  = & \langle - iky (M_{k}f), g \rangle_{L^{2}_{y}} + \langle M_{k}(iky f), g \rangle_{L^{2}_{y}} \\
  = & \langle k\partial_{\xi}M_{k}\widehat{f} + kM_{k}\partial_{\xi}\widehat{f},\widehat{g} \rangle_{L^{2}_{\xi}} + \langle -k M_{k}\partial_{\xi}\widehat{f}, \widehat{g} \rangle_{L^{2}_{\xi}} \\
  = & \langle (k\partial_{\xi}M_{k})(D_{y})f ,g \rangle_{L^{2}_{y}}\\
  =&  \langle \eta^{\frac{1}{3}}\vert k \vert^{\frac{2}{3}} \varphi^{'}(\eta^{\frac{1}{3}} \vert k \vert^{-\frac{1}{3}} {\rm sgn}(k) D_{y})f, g \rangle_{L^{2}_{y}},
\end{align*}
we obtain
\begin{align*}
  \frac{d}{dt} \langle M_{k}\theta_{k}, \theta_{k} \rangle_{L^{2}_{y}} + \langle 2\eta(D_{y}^{2} +  k^{2})\varphi(\eta^{\frac{1}{3}}\vert k \vert^{-\frac{1}{3}}{\rm sgn}(k)D_{y})\theta_{k}, \theta_{k} \rangle_{L^{2}_{y}} \\
  + \langle \eta^{\frac{1}{3}}\vert k \vert^{\frac{2}{3}}\varphi^{'}(\eta^{\frac{1}{3}}\vert k \vert^{-\frac{1}{3}}{\rm sgn}(k)D_{y})\theta_{k}, \theta_{k} \rangle_{L^{2}_{y}} = 0.
\end{align*}
Together with (\ref{linear_equation_identity1}), this gives
\begin{equation}\label{linear_equation_identity1_together}
\begin{split}
  \frac{d}{dt} \left(\langle M_{k}\theta_{k}, \theta_{k} \rangle_{L^{2}_{y}} + \Vert \theta_{k} \Vert_{L_{y}^{2}}^{2}\right)   + \left\langle \left(2\eta(D_{y}^{2} +  k^{2})\big(\varphi(\eta^{\frac{1}{3}}\vert k \vert^{-\frac{1}{3}}{\rm sgn}(k)D_{y})+1\big) \right.\right.\\
  \left.\left. \qquad\qquad  + \eta^{\frac{1}{3}}\vert k \vert^{\frac{2}{3}}\varphi^{'}(\eta^{\frac{1}{3}}\vert k \vert^{-\frac{1}{3}}{\rm sgn}(k)D_{y})\right) \theta_{k}, \theta_{k} \right\rangle_{L^{2}_{y}} = 0.
\end{split}
\end{equation}
By the choice of the function $\varphi$, there holds
\[
  \eta(\xi^{2} +  k^{2})\big(1+2\varphi(\eta^{\frac{1}{3}}\vert k \vert^{-\frac{1}{3}}{\rm sgn}(k)\xi)\big) + \eta^{\frac{1}{3}} \vert k \vert^{\frac{2}{3}}\varphi^{'}(\eta^{\frac{1}{3}}\vert k \vert^{-\frac{1}{3}}{\rm sgn}(k)\xi) \geq \frac{1}{4} \eta^{\frac{1}{3}}\vert k \vert^{\frac{2}{3}},
\]
for all $k \in \mathbb{Z}, \ k\neq 0,  \eta > 0,  \xi \in \mathbb{R}$. In fact, when $\vert \eta^{\frac{1}{3}}\vert k \vert^{-\frac{1}{3}}{\rm sgn}(k)\xi \vert \leq 1$, we have \[ 
\varphi^{'}(\eta^{\frac{1}{3}}\vert k \vert^{-\frac{1}{3}}{\rm sgn}(k)\xi) = \frac{1}{4},
\]
and the above inequality clearly holds. When $\vert \eta^{\frac{1}{3}}\vert k \vert^{-\frac{1}{3}}{\rm sgn}(k)\xi \vert > 1$, we have 
\[
\frac{1}{4} \eta^{\frac{1}{3}}\vert k \vert^{\frac{2}{3}} \leq \frac{1}{4} \eta^{\frac{1}{3}}\vert k \vert^{\frac{2}{3}} \vert \eta^{\frac{1}{3}}\vert k \vert^{-\frac{1}{3}}{\rm sgn}(k)\xi \vert^{2} \leq  \frac{1}{4}\eta \xi^{2}.
\]
 According to
\begin{align*}
  & \left\langle \left(2\eta(D_{y}^{2} +  k^{2})\big(\varphi(\eta^{\frac{1}{3}}\vert k \vert^{-\frac{1}{3}}{\rm sgn}(k)D_{y})+1\big) + \eta^{\frac{1}{3}}\vert k \vert^{\frac{2}{3}}\varphi^{'}(\eta^{\frac{1}{3}}\vert k \vert^{-\frac{1}{3}}{\rm sgn}(k)D_{y})\right) \theta_{k}, \theta_{k} \right\rangle_{L^{2}_{y}} \\
  & = \left\langle \left(\eta(\xi^{2} +  k^{2})\big(2\varphi(\eta^{\frac{1}{3}}\vert k \vert^{-\frac{1}{3}}{\rm sgn}(k)\xi)+1\big) + \eta^{\frac{1}{3}}\vert k \vert^{\frac{2}{3}}\varphi^{'}(\eta^{\frac{1}{3}}\vert k \vert^{-\frac{1}{3}}{\rm sgn}(k)\xi)\right) \widehat{\theta}_{k}, \widehat{\theta_{k}} \right\rangle_{L^{2}_{\xi}}  \\
  & \qquad + \big\langle \eta(\xi^{2} +  k^{2}) \widehat{\theta_{k}}, \widehat{\theta_{k}} \big\rangle_{L^{2}_{\xi}} \\
  & \geq \big\langle \frac{1}{4}\eta^{\frac{1}{3}}\vert k \vert^{\frac{2}{3}} \widehat{\theta_{k}}, \widehat{\theta_{k}} \big\rangle_{L^{2}_{\xi}} + \big\langle \eta(\xi^{2} +  k^{2}) \widehat{\theta_{k}}, \widehat{\theta_{k}} \big\rangle_{L^{2}_{\xi}} \\
  & = \frac{1}{4}\eta^{\frac{1}{3}}\vert k \vert^{\frac{2}{3}} \Vert \theta_{k} \Vert^{2}_{L^{2}_{y}} + \eta \Vert D_{y}\theta_{k} \Vert^{2}_{L^{2}_{y}} + \eta  k^{2} \Vert \theta_{k} \Vert^{2}_{L^{2}_{y}},
\end{align*}
together with (\ref{linear_equation_identity1_together}), we get
\begin{equation*}\label{theta_k_inequality}
  \frac{d}{dt} \left(\langle M_{k}\theta_{k}, \theta_{k} \rangle_{L^{2}_{y}} + \Vert \theta_{k} \Vert_{L_{y}^{2}}^{2}\right) + \frac{1}{4}\eta^{\frac{1}{3}}\vert k \vert^{\frac{2}{3}} \Vert \theta_{k} \Vert^{2}_{L^{2}_{y}} + \eta \Vert D_{y}\theta_{k} \Vert^{2}_{L^{2}_{y}} + \eta  k^{2} \Vert \theta_{k} \Vert^{2}_{L^{2}_{y}} \leq 0.
\end{equation*}
We see
\begin{equation}\label{theta_inequation}
\begin{split}
 &\frac{d}{dt} \left(\Vert \sqrt{M_{k} + 1}\widehat{\theta_{k}}(t) \Vert^{2}_{L^{2}_{\xi}} e^{\frac{1}{8}\eta^{\frac{1}{3}}\vert k \vert^{\frac{2}{3}}t}\right) \\&=  e^{\frac{1}{8}\eta^{\frac{1}{3}}\vert k \vert^{\frac{2}{3}}t}\left( \frac{d}{dt} \left(\Vert \sqrt{M_{k} + 1}\widehat{\theta_{k}}(t) \Vert^{2}_{L^{2}_{\xi}}\right) + \frac{1}{8}\eta^{\frac{1}{3}}\vert k \vert^{\frac{2}{3}} \Vert \sqrt{M_{k} + 1}\widehat{\theta_{k}}(t) \Vert^{2}_{L^{2}_{\xi}} \right) \\
 & \leq  e^{\frac{1}{8}\eta^{\frac{1}{3}}\vert k \vert^{\frac{2}{3}}t}\left( \frac{d}{dt} \left( \langle M_{k}\theta_{k}, \theta_{k} \rangle_{L^{2}_{y}} + \Vert \theta_{k} \Vert_{L_{y}^{2}}^{2} \right) + \frac{1}{4}\eta^{\frac{1}{3}}\vert k \vert^{\frac{2}{3}} \Vert \theta_{k} \Vert^{2}_{L^{2}_{y}} \right) \\
  &\leq  \, 0.
\end{split}
\end{equation}
Integrating   (\ref{theta_inequation}) in $t$ and using the properties of $M_{k}$, we obtain
\begin{align}\label{theta}
 \Vert \theta_{k}(t) \Vert_{L_{y}^{2}} \leq \sqrt{2} \Vert \theta_{k}(0) \Vert_{L^{2}_{y}}e^{-\frac{1}{16}\eta^{\frac{1}{3}}\vert k \vert^{\frac{2}{3}}t}. 
\end{align}
When $k = 0$, for (\ref{linear_equation_identity1}), we get   $\frac{1}{2}\frac{d}{dt}\Vert \theta_{0}(t) \Vert^{2}_{L^{2}_{y}} \leq 0 $, which implies that $\Vert \theta_{0}(t) \Vert_{L_{y}^{2}} \leq \Vert \theta_{0}(0) \Vert_{L^{2}_{y}}$.

Differentiating the third equation in (\ref{linear_system__projecting_y}) with respect to $y$ leads to
\begin{align*}
  \partial_{t}D_{y}^{N}\theta_{k} + iky D_{y}^{N}\theta_{k} + kN D_{y}^{N-1}\theta_{k} + \eta (D^{2}_{y}+ k^{2})D_{y}^{N} \theta_{k} = 0.
\end{align*}
Taking the $L_{y}^{2}$-inner product with $(1+M_{k})D_{y}^{N}\theta_{k}$ gives
\begin{align*}
  & \frac{d}{dt} \left(\langle M_{k}D_{y}^{N}\theta_{k}, D_{y}^{N}\theta_{k} \rangle_{L^{2}_{y}} + \Vert D_{y}^{N}\theta_{k} \Vert_{L_{y}^{2}}^{2}\right) + \frac{1}{4}\eta^{\frac{1}{3}}\vert k \vert^{\frac{2}{3}} \Vert D_{y}^{N}\theta_{k} \Vert^{2}_{L^{2}_{y}} \\
  &\qquad + \eta \Vert D^{N+1}_{y}\theta_{k} \Vert^{2}_{L^{2}_{y}} + \eta  k^{2} \Vert D_{y}^{N}\theta_{k} \Vert^{2}_{L^{2}_{y}} \\
  & \leq -2{\rm  Re} \langle kND_{y}^{N-1}\theta_{k}, (1+M_{k})D_{y}^{N}\theta_{k} \rangle_{L^{2}_{y}} \\
  & \leq \frac{1}{8}\eta^{\frac{1}{3}} \vert k \vert^{\frac{2}{3}} \Vert D_{y}^{N}\theta_{k} \Vert_{L^{2}_{y}}^{2} + 32N^{2}\eta^{-\frac{1}{3}} \vert k \vert^{\frac{4}{3}} \Vert D_{y}^{N-1}\theta_{k} \Vert_{L^{2}_{y}}^{2},
  \end{align*}
  therefore,
  \begin{align*}
 & \frac{d}{dt} \left(\langle M_{k}D_{y}^{N}\theta_{k}, D_{y}^{N}\theta_{k} \rangle_{L^{2}_{y}} + \Vert D_{y}^{N}\theta_{k} \Vert_{L_{y}^{2}}^{2}\right) + \frac{1}{8}\eta^{\frac{1}{3}}\vert k \vert^{\frac{2}{3}} \Vert D_{y}^{N}\theta_{k} \Vert^{2}_{L^{2}_{y}} \\
 & \qquad+ \eta \Vert D^{N+1}_{y}\theta_{k} \Vert^{2}_{L^{2}_{y}} + \eta  k^{2} \Vert D_{y}^{N}\theta_{k} \Vert^{2}_{L^{2}_{y}} \\
 & \leq 32N^{2}\eta^{-\frac{1}{3}} \vert k \vert^{\frac{4}{3}} \Vert D_{y}^{N-1}\theta_{k} \Vert_{L^{2}_{y}}^{2}.
\end{align*}
Similarly, we get
\begin{align*}
  \frac{d}{dt} \left(\Vert \sqrt{M_{k} + 1}\xi^{N}\widehat{\theta_{k}}(t) \Vert^{2}_{L^{2}_{\xi}} e^{\frac{1}{16}\eta^{\frac{1}{3}}\vert k \vert^{\frac{2}{3}}t}\right) 
  \leq e^{\frac{1}{16}\eta^{\frac{1}{3}}\vert k \vert^{\frac{2}{3}}t} 32N^{2}\eta^{-\frac{1}{3}} \vert k \vert^{\frac{4}{3}} \Vert D_{y}^{N-1}\theta_{k} \Vert_{L^{2}_{y}}^{2},
\end{align*}
which leads to 
\begin{align*}
   \Vert D_{y}^{N}\theta_{k}(t) \Vert_{L_{y}^{2}}^{2}  \leq e^{-\frac{1}{16}\eta^{\frac{1}{3}}\vert k \vert^{\frac{2}{3}}t}
   \Big(2 \Vert D_{y}^{N}\theta_{k}(0) \Vert^{2}_{L^{2}_{y}}+  32N^{2} \int_{0}^{t} \eta^{-\frac{1}{3}} \vert k \vert^{\frac{4}{3}} \Vert D_{y}^{N-1}\theta_{k}(s) \Vert_{L^{2}_{y}}^{2} e^{\frac{1}{16}\eta^{\frac{1}{3}}\vert k \vert^{\frac{2}{3}}s} ds\Big).
\end{align*}
Now we prove the  inequality (\ref{the_linear_stability_result2}) by induction. For $N = 1$,
\begin{align*}
  \Vert D_{y}\theta_{k}(t) \Vert_{L_{y}^{2}}^{2} &\leq  e^{-\frac{1}{16}\eta^{\frac{1}{3}}\vert k \vert^{\frac{2}{3}}t} \Big(  2 \Vert D_{y}\theta_{k}(0) \Vert^{2}_{L^{2}_{y}} +  32\int_{0}^{t} \eta^{-\frac{1}{3}} \vert k \vert^{\frac{4}{3}} \Vert \theta_{k}(s) \Vert_{L^{2}_{y}}^{2} e^{\frac{1}{16}\eta^{\frac{1}{3}}\vert k \vert^{\frac{2}{3}}s} ds \Big) \\
  &  \leq  e^{-\frac{1}{16}\eta^{\frac{1}{3}}\vert k \vert^{\frac{2}{3}}t} \Big(  2 \Vert D_{y}\theta_{k}(0) \Vert^{2}_{L^{2}_{y}}
   +  64
    \eta^{-\frac{1}{3}} \vert k \vert^{\frac{4}{3}} 
     \Vert \theta_{k}(0) \Vert_{L^{2}_{y}}^{2}
    \int_{0}^{t} e^{-\frac{1}{16}\eta^{\frac{1}{3}}\vert k \vert^{\frac{2}{3}}s} ds \Big) \\
  & \leq e^{-\frac{1}{16}\eta^{\frac{1}{3}}\vert k \vert^{\frac{2}{3}}t}  \left(  2 \Vert D_{y}\theta_{k}(0) \Vert^{2}_{L^{2}_{y}} +  1024 (\eta^{-1}\vert k \vert)^{\frac{2}{3}}\Vert \theta_{k}(0) \Vert_{L^{2}_{y}}^{2} \right),
  \end{align*}
  which leads to
  \begin{align*}
 \Vert D_{y}\theta_{k}(t) \Vert_{L_{y}^{2}} \leq C_{1} e^{-\frac{1}{32}\eta^{\frac{1}{3}}\vert k \vert^{\frac{2}{3}}t} \left( \Vert D_{y}\theta_{k}(0) \Vert_{L^{2}_{y}} + (\eta^{-1}\vert k \vert)^{\frac{1}{3}}\Vert \theta_{k}(0) \Vert_{L^{2}_{y}} \right).
\end{align*}
Assume that for $n \leq N$, there exist two positive constants $C_{n}>0$ and $c_{n}>0$, such that
\[
  \Vert D_{y}^{n}\theta_{k}(t) \Vert_{L_{y}^{2}} \leq C_{n} e^{-c_{n}\eta^{\frac{1}{3}}\vert k \vert^{\frac{2}{3}}t} \left( \Vert D^{n}_{y}\theta_{k}(0) \Vert_{L^{2}_{y}} + (\eta^{-1}\vert k \vert)^{\frac{n}{3}}\Vert \theta_{k}(0) \Vert_{L^{2}_{y}} \right).
\]
Then for $n=N+1$, we have
\begin{align*}
 \Vert D_{y}^{N+1}\theta_{k}(t) \Vert_{L_{y}^{2}}^{2}
 &  \leq  e^{-\frac{1}{16}\eta^{\frac{1}{3}}\vert k \vert^{\frac{2}{3}}t}\Big( 2 \Vert D_{y}^{N+1}\theta_{k}(0) \Vert^{2}_{L^{2}_{y}} \\
 &  \qquad +  32(N+1)^{2} \int_{0}^{t} \eta^{-\frac{1}{3}} \vert k \vert^{\frac{4}{3}} \Vert D_{y}^{N}\theta_{k}(s) \Vert_{L^{2}_{y}}^{2} e^{\frac{1}{16}\eta^{\frac{1}{3}}\vert k \vert^{\frac{2}{3}}s} ds \Big).
 \end{align*}
When $\frac{1}{16}-2c_{N} < 0$, we have
\begin{align*}
  & \int_{0}^{t} \eta^{-\frac{1}{3}} \vert k \vert^{\frac{4}{3}} \Vert D_{y}^{N}\theta_{k}(s) \Vert_{L^{2}_{y}}^{2} e^{\frac{1}{16}\eta^{\frac{1}{3}}\vert k \vert^{\frac{2}{3}}s} ds \\
  &  \leq C_{N}^{2} \eta^{-\frac{1}{3}} \vert k \vert^{\frac{4}{3}} \left( \Vert D^{N}_{y}\theta_{k}(0) \Vert_{L^{2}_{y}} + (\eta^{-1}\vert k \vert)^{\frac{N}{3}}\Vert \theta_{k}(0) \Vert_{L^{2}_{y}} \right)^{2} \int_{0}^{t} e^{(\frac{1}{16}-2c_{N})\eta^{\frac{1}{3}}\vert k \vert^{\frac{2}{3}}s} ds \\
  &  \leq \frac{2C_{N}^{2}}{2c_{N}-\frac{1}{16}}(\eta^{-1}\vert k \vert)^{\frac{2}{3}} \left( \Vert D^{N}_{y}\theta_{k}(0) \Vert^{2}_{L^{2}_{y}} + (\eta^{-1}\vert k \vert)^{\frac{2N}{3}}\Vert \theta_{k}(0) \Vert^{2}_{L^{2}_{y}} \right).
\end{align*}
Because
\begin{align*}
  (\eta^{-1}\vert k \vert)^{\frac{2}{3}}\Vert D^{N}_{y}\theta_{k}(0) \Vert^{2}_{L^{2}_{y}} & \leq (\eta^{-1}\vert k \vert)^{\frac{2}{3}} \Vert \theta_{k}(0) \Vert^{\frac{2}{N+1}}_{L^{2}_{y}} \Vert D_{y}^{N+1}\theta_{k}(0) \Vert^{\frac{2N}{N+1}}_{L^{2}_{y}} \\
  & \leq \frac{1}{N+1} \left( (\eta^{-1}\vert k \vert)^{\frac{2N+2}{3}} \Vert \theta_{k}(0) \Vert^{2}_{L^{2}_{y}} \right) + \frac{N}{N+1} \Vert D_{y}^{N+1}\theta_{k}(0) \Vert^{2}_{L^{2}_{y}},
\end{align*}
 we have
\begin{equation}\label{1_2cn_leq_0}
   \Vert D_{y}^{N+1}\theta_{k}(t) \Vert_{L_{y}^{2}} \leq 
    C_{N+1} e^{-\frac{1}{32}\eta^{\frac{1}{3}}\vert k \vert^{\frac{2}{3}}t} 
   \left( \Vert D_{y}^{N+1}\theta_{k}(0) \Vert_{L^{2}_{y}} + (\eta^{-1}\vert k \vert)^{\frac{N+1}{3}} \Vert \theta_{k}(0) \Vert_{L^{2}_{y}} \right).
\end{equation}
When $0 \leq \frac{1}{16}-2c_{N} < \alpha_{N} < \frac{1}{16}$, we get
\begin{align*}
  & \int_{0}^{t} \eta^{-\frac{1}{3}} \vert k \vert^{\frac{4}{3}} \Vert D_{y}^{N}\theta_{k}(s) \Vert_{L^{2}_{y}}^{2} e^{\frac{1}{16}\eta^{\frac{1}{3}}\vert k \vert^{\frac{2}{3}}s} ds \\
  & \leq C_{N}^{2} \eta^{-\frac{1}{3}} \vert k \vert^{\frac{4}{3}} \left( \Vert D^{N}_{y}\theta_{k}(0) \Vert_{L^{2}_{y}} + (\eta^{-1}\vert k \vert)^{\frac{N}{3}}\Vert \theta_{k}(0) \Vert_{L^{2}_{y}} \right)^{2} \int_{0}^{t} e^{\alpha_{N}\eta^{\frac{1}{3}}\vert k \vert^{\frac{2}{3}}s} ds \\
  &  \leq 2 C_{N}^{2}\eta^{-\frac{1}{3}} \vert k \vert^{\frac{4}{3}} \left( \Vert D^{N}_{y}\theta_{k}(0) \Vert^{2}_{L^{2}_{y}} + (\eta^{-1}\vert k \vert)^{\frac{2N}{3}}\Vert \theta_{k}(0) \Vert^{2}_{L^{2}_{y}} \right)\frac{e^{\alpha_{N}\eta^{\frac{1}{3}}\vert k \vert^{\frac{2}{3}}t}}{\alpha_{N}\eta^{\frac{1}{3}}\vert k \vert^{\frac{2}{3}}} \\
  & \leq C_{N+1}^{2} \left( (\eta^{-1}\vert k \vert)^{\frac{2}{3}}\Vert D^{N}_{y}\theta_{k}(0) \Vert^{2}_{L^{2}_{y}} + (\eta^{-1}\vert k \vert)^{\frac{2N+2}{3}}\Vert \theta_{k}(0) \Vert^{2}_{L^{2}_{y}} \right)e^{\alpha_{N}\eta^{\frac{1}{3}}\vert k \vert^{\frac{2}{3}}t} \\
  &  \leq C_{N+1}^{2} (\Vert D_{y}^{N+1}\theta_{k}(0) \Vert_{L^{2}_{y}}^{2} + (\eta^{-1}\vert k \vert)^{\frac{2N+2}{3}} \Vert \theta_{k}(0) \Vert^{2}_{L^{2}_{y}})e^{\alpha_{N}\eta^{\frac{1}{3}}\vert k \vert^{\frac{2}{3}}t},
\end{align*}
where $\alpha_{N}$ is a positive constant and $C_{N+1}$ is  a generic positive constant which may change from line to line, and we will often do so without any remark.
Hence,
\begin{align*}
   \Vert D_{y}^{N+1}\theta_{k}(t) \Vert_{L_{y}^{2}}^{2}
    &
   \leq  e^{-\frac{1}{16}\eta^{\frac{1}{3}}\vert k \vert^{\frac{2}{3}}t}\Big( 2 \Vert D_{y}^{N+1}\theta_{k}(0) \Vert^{2}_{L^{2}_{y}} \\
   & \qquad +  32(N+1)^{2} \int_{0}^{t} \eta^{-\frac{1}{3}} \vert k \vert^{\frac{4}{3}} \Vert D_{y}^{N}\theta_{k}(s) \Vert_{L^{2}_{y}}^{2} e^{\frac{1}{16}\eta^{\frac{1}{3}}\vert k \vert^{\frac{2}{3}}s} ds \Big) \\
  &  \leq 2e^{-\frac{1}{16}\eta^{\frac{1}{3}}\vert k \vert^{\frac{2}{3}}t}\Vert D_{y}^{N+1}\theta_{k}(0) \Vert^{2}_{L^{2}_{y}}
   + 32(N+1)^{2} C^{2}_{N+1} \Big(\Vert D_{y}^{N+1}\theta_{k}(0) \Vert_{L^{2}_{y}}^{2} \\
  & \qquad + (\eta^{-1}\vert k \vert)^{\frac{2N+2}{3}} \Vert \theta_{k}(0) \Vert^{2}_{L^{2}_{y}}\Big)
  e^{(\alpha_{N}-\frac{1}{16})\eta^{\frac{1}{3}}\vert k \vert^{\frac{2}{3}}t} \\
  & \leq C_{N+1}^{2}e^{(\alpha_{N}-\frac{1}{16})\eta^{\frac{1}{3}}\vert k \vert^{\frac{2}{3}}t}\left( \Vert D_{y}^{N+1}\theta_{k}(0) \Vert_{L^{2}_{y}}^{2} + (\eta^{-1}\vert k \vert)^{\frac{2N+2}{3}} \Vert \theta_{k}(0) \Vert^{2}_{L^{2}_{y}} \right).
\end{align*}
By the estimate  (\ref{1_2cn_leq_0}) and letting $c_{N+1} = \frac{1}{32} - \frac{1}{2}\alpha_{N} $, we can show 
\begin{equation*}
   \Vert D_{y}^{N+1}\theta_{k}(t) \Vert_{L_{y}^{2}} \leq C_{N+1}e^{-c_{N+1}\eta^{\frac{1}{3}}\vert k \vert^{\frac{2}{3}}t}\left( \Vert D_{y}^{N+1}\theta_{k}(0) \Vert_{L^{2}_{y}} + (\eta^{-1}\vert k \vert)^{\frac{N+1}{3}} \Vert \theta_{k}(0) \Vert_{L^{2}_{y}} \right).
\end{equation*}

Now, we prove the inequality (\ref{the_linear_stability_result1}). For $w_{k}$ and  $j_{k}$, we have
\begin{align*}
& \frac{d}{dt} \Vert w_{k} \Vert^{2}_{L_{y}^{2}} + \langle 2 \nu (D_{y}^{2} +  k^{2}) w_{k}, w_{k} \rangle_{L_{y}^{2}}= 2 {\rm  Re} \langle ik \theta_{k}, w_{k} \rangle_{L_{y}^{2}} + 2 {\rm  Re} \langle ik j_{k}, w_{k} \rangle_{L_{y}^{2}}, 
\\
  & \frac{d}{dt} \Vert j_{k} \Vert^{2}_{L_{y}^{2}} + \langle 2 \mu (D_{y}^{2} +  k^{2}) j_{k}, j_{k} \rangle_{L_{y}^{2}} = 2 {\rm  Re} \langle ik w_{k}, j_{k} \rangle_{L_{y}^{2}} + 2 {\rm  Re} \langle 2 ik \partial_{y}(k^{2}+D_{y}^{2})^{-1} j_{k}, j_{k} \rangle_{L_{y}^{2}}.
\end{align*}
And we define function $\phi_{k}$ for $k \neq 0$ as follows, 
\begin{equation*}
\phi_{k}(\xi) = 
\begin{cases}
  \frac{6(k^{2} + \xi_{0}^{2})^{2}}{k^{4}} - (2 + \pi) ,\qquad & \xi > 0, \\
  \frac{6(k^{2} + \xi_{0}^{2})^{2}}{(k^{2} + \xi^{2})^{2}} - (2 + \pi), \qquad & \xi \in [-\xi_{0},0], \\
  (4 - \pi) e^{\frac{24\xi_{0}(\xi+\xi_{0})}{(4-\pi)(k^{2}+\xi_{0}^{2})}}, \qquad & \xi \in (-\infty, -\xi_{0}),
\end{cases}
\end{equation*}
where $\xi_{0}$ is a real positive solution of the equation $\nu \xi_{0} (k^{2} + \xi_{0}^{2}) = 96 \vert k \vert$. 
It is easy to check that $\phi_{k}(\xi) \in C^{1}(\mathbb{R})$ for  $k \neq 0$. Therefore, we get $\xi_{0} \leq 96 \nu^{-1}\vert k \vert^{-1} \leq 96 \nu^{-1}$, which implies that  $0 < \phi_{k} (\xi)\leq C\nu^{-4}$ and $0 \leq \phi_{k}^{'}(\xi) \leq C\frac{\nu^{-3}}{\vert k \vert}$, $\xi\in \mathbb{R}$, where $C$ is a positive constant.
And  we define
\[ M^{'}_{k} = \varphi(\nu^{\frac{1}{3}}\vert k \vert^{-\frac{1}{3}}{\rm sgn}(k)D_{y}) + \phi_{k}({\rm sgn}(k)D_{y}), ~{\rm for}\, k \neq 0, \qquad M^{'}_{0}= 0.
\]
Then we can make a conclusion
 \begin{align*}
 &  \nu (\xi^{2} +  k^{2})(1+2M^{'}_{k}(\xi)) + k\partial_{\xi}M^{'}_{k}(\xi)  \geq \frac{1}{4}\nu^{\frac{1}{3}}\vert k \vert^{\frac{2}{3}},
\\
 &  \nu (\xi^{2} +  k^{2})(1+2M^{'}_{k}(\xi)) + k\partial_{\xi}M^{'}_{k}(\xi) + (1+M^{'}_{k}(\xi)) \frac{4k\xi}{k^{2}+\xi^{2}} \geq \frac{1}{4}\nu^{\frac{1}{3}}\vert k \vert^{\frac{2}{3}},
   \end{align*}
 with $ k\in \mathbb{Z}$ and $k \neq 0$.
 Multiplying the $w_{k}$ equation by $(1+M^{'}_{k})w_{k}$ and the $j_{k}$ equation by $(1+M^{'}_{k})j_{k}$, we can get
\begin{equation*}\label{wk_jk_m_inequation}
\begin{split}
  & \frac{d}{dt}\langle (1+M^{'}_{k})w_{k},w_{k} \rangle_{L_{y}^{2}} + \nu \Vert D_{y}w_{k} \Vert^{2}_{L_{y}^{2}} +  \nu k^{2} \Vert w_{k} \Vert^{2}_{L_{y}^{2}} + \frac{1}{4}\nu^{\frac{1}{3}}\vert k \vert^{\frac{2}{3}}\Vert w_{k} \Vert^{2}_{L_{y}^{2}} \\
  & \quad+ \frac{d}{dt}\langle (1+M^{'}_{k})j_{k},j_{k} \rangle_{L_{y}^{2}} + \nu \Vert D_{y}j_{k} \Vert^{2}_{L_{y}^{2}} +  \nu k^{2} \Vert j_{k} \Vert^{2}_{L_{y}^{2}} + \frac{1}{4}\nu^{\frac{1}{3}}\vert k \vert^{\frac{2}{3}}\Vert j_{k} \Vert^{2}_{L_{y}^{2}} \\
 & \leq 2 {\rm  Re} \langle ik\theta_{k}, (1+M^{'}_{k})w_{k} \rangle_{L_{y}^{2}}.
\end{split}
\end{equation*}
Applying Young's inequality to the right-hand side and removing some terms on the left-hand side yields
\begin{align*}
  & \frac{d}{dt}\left(\langle (1+M^{'}_{k})w_{k},w_{k} \rangle_{L_{y}^{2}} + \langle (1+M^{'}_{k})j_{k},j_{k} \rangle_{L_{y}^{2}}\right) + \frac{1}{8}\nu^{\frac{1}{3}}\vert k \vert^{\frac{2}{3}}\Vert\big( w_{k} ,j_{k}\big) \Vert^{2}_{L_{y}^{2}} \\
 & \leq C\nu^{-\frac{1}{3} - 8} \vert k \vert^{\frac{4}{3}}\Vert \theta_{k} \Vert^{2}_{L_{y}^{2}},
\end{align*}
which gives that 
\begin{align*}
  & \frac{d}{dt} \left(\Big(\Vert \sqrt{1+M^{'}_{k}} w_{k} \Vert_{L_{y}^{2}}^{2} + \Vert \sqrt{1+M^{'}_{k}} j_{k} \Vert_{L_{y}^{2}}^{2}\Big )e^{\frac{1}{8(1+M_{k}^{'})} \nu^{\frac{1}{3}}\vert k \vert^{\frac{2}{3}}t}\right)  \\
   &\leq  e^{\frac{1}{8(1+M^{'}_{k})} \nu^{\frac{1}{3}}\vert k \vert^{\frac{2}{3}}t} \left( \frac{d}{dt}\Big(\Vert \sqrt{1+M^{'}_{k}} w_{k} \Vert_{L_{y}^{2}}^{2} + \Vert \sqrt{1+M^{'}_{k}} j_{k} \Vert_{L_{y}^{2}}^{2} \Big)
   + \frac{1}{8} \nu^{\frac{1}{3}}\vert k \vert^{\frac{2}{3}}\Vert \big(w_{k}, j_{k}\big) \Vert_{L_{y}^{2}}^{2} \right) \\  
   & \leq Ce^{\frac{1}{8(1+M^{'}_{k})} \nu^{\frac{1}{3}}\vert k \vert^{\frac{2}{3}}t} \nu^{-\frac{1}{3}-8} \vert k \vert^{\frac{4}{3}}\Vert \theta_{k} \Vert^{2}_{L_{y}^{2}} .
   \end{align*}
Integrating in $t$  and using the  inequality (\ref{theta}), we obtain
  \begin{align*}
 \Vert\big( w_{k}(t), j_{k}(t)\big) \Vert^{2}_{L_{y}^{2}} 
 \leq C^{2} \left( \nu^{-4}\Vert\big (w_{k}(0),j_{k}(0) \big) \Vert^{2}_{L_{y}^{2}}  + \nu^{-12}(\nu^{-1} \vert k \vert)^{\frac{2}{3}}\Vert \theta_{k}(0) \Vert^{2}_{L_{y}^{2}} \right)e^{- 2c\nu^{\frac{13}{3}}\vert k \vert^{\frac{2}{3}}t}.
\end{align*}

Differentiating the $w_{k}$ equation and the  $j_{k}$ equation in (\ref{linear_system__projecting_y}) with respect to $y$, respectively, leads to
\begin{equation*}
  \partial_{t}D_{y}^{N}w_{k} + iky D_{y}^{N}w_{k} + kN D_{y}^{N-1}w_{k} + \nu (D^{2}_{y}+ k^{2})D_{y}^{N} w_{k} = ik D_{y}^{N}\theta_{k} + ik D_{y}^{N}j_{k},
\end{equation*}
\begin{equation*}
 \partial_{t}D_{y}^{N}j_{k} + iky D_{y}^{N}j_{k} + kN D_{y}^{N-1}j_{k} + \mu (D^{2}_{y}+ k^{2})D_{y}^{N} j_{k} -2ik D_{y}^{N}b_{k}^{1} = ik D_{y}^{N}w_{k}.
\end{equation*}
Taking the $L^{2}_{y}$-inner product with $(1+M^{'}_{k})D^{N}_{y}w_{k}$ and $(1+M^{'}_{k})D^{N}_{y}j_{k}$, we get
\begin{align*}
	& \frac{d}{dt}\left(\langle (1+M^{'}_{k})D^{N}_{y}w_{k}, D^{N}_{y}w_{k} \rangle_{L_{y}^{2}} + \langle (1+M^{'}_{k})D^{N}_{y}j_{k}, D^{N}_{y}j_{k} \rangle_{L_{y}^{2}}\right) \\
	& \qquad + \frac{1}{4}\nu^{\frac{1}{3}}\vert k \vert^{\frac{2}{3}}
	 \Vert\big( D^{N}_{y}w_{k} , D^{N}_{y}j_{k}\big) \Vert_{L_{y}^{2}}^{2} \\
	& \leq -2 {\rm  Re} \langle kN D_{y}^{N-1}w_{k}, (1+M^{'}_{k})D_{y}^{N}w_{k}\rangle_{L_{y}^{2}} - 2 {\rm  Re} \langle kN D_{y}^{N-1}j_{k}, (1+M^{'}_{k})D_{y}^{N}j_{k}\rangle_{L_{y}^{2}} \\
	& \qquad + 2 {\rm  Re} \langle ik D_{y}^{N}\theta_{k}, (1+M^{'}_{k})D_{y}^{N}w_{k}\rangle_{L_{y}^{2}} \\
	& \leq \frac{1}{16}\nu^{\frac{1}{3}}\vert k \vert^{\frac{2}{3}} \Vert \big(D^{N}_{y}w_{k} , D^{N}_{y}j_{k} \big)\Vert_{L_{y}^{2}}^{2} 
	+ C N^{2} \nu^{-\frac{1}{3}-8}\vert k \vert^{\frac{4}{3}}\Vert\big( D^{N-1}_{y}w_{k} , D ^{N-1}_{y}j_{k} \big)\Vert_{L_{y}^{2}}^{2} \\
	& \qquad + \frac{1}{16}\nu^{\frac{1}{3}}\vert k \vert^{\frac{2}{3}}\Vert D^{N}_{y}w_{k} \Vert_{L_{y}^{2}}^{2} + C \nu^{-\frac{1}{3}-8}\vert k \vert^{\frac{4}{3}}\Vert D^{N}_{y}\theta_{k} \Vert_{L_{y}^{2}}^{2},
\end{align*}
which gives that
\begin{align*}
& \frac{d}{dt}\left(\langle (1+M^{'}_{k})D^{N}_{y}w_{k}, D^{N}_{y}w_{k} \rangle_{L_{y}^{2}} + \langle (1+M^{'}_{k})D^{N}_{y}j_{k}, D^{N}_{y}j_{k} \rangle_{L_{y}^{2}}\right) \\
 & \qquad + \frac{1}{8}\nu^{\frac{1}{3}}\vert k \vert^{\frac{2}{3}} \Vert\big( D^{N}_{y}w_{k} , D^{N}_{y}j_{k} \big)\Vert_{L_{y}^{2}}^{2} \\
	& \leq C N^{2} \nu^{-\frac{1}{3}-8}\vert k \vert^{\frac{4}{3}}\Vert\big( D^{N-1}_{y}w_{k} , D^{N-1}_{y}j_{k} \big)\Vert_{L_{y}^{2}}^{2}+ C \nu^{-\frac{1}{3}-8}\vert k \vert^{\frac{4}{3}}\Vert D^{N}_{y}\theta_{k} \Vert_{L_{y}^{2}}^{2}.
\end{align*}
Then we have
\begin{align*}
 &  \Vert D^{N}_{y}w_{k}(t) \Vert_{L_{y}^{2}}^{2} 
 +\Vert D^{N}_{y}j_{k}(t) \Vert_{L_{y}^{2}}^{2} \\
  & \leq  C\nu^{-4}\Vert \big(D^{N}_{y}w_{k}(0) , D^{N}_{y}j_{k}(0)\big) \Vert_{L_{y}^{2}}^{2}e^{-\frac{1}{8(1+M^{'}_{k})}\nu^{\frac{1}{3}}\vert k \vert^{\frac{2}{3}}t} \\
  & \qquad +  \int_{0}^{t} C N^{2} \nu^{-\frac{1}{3}-8}\vert k \vert^{\frac{4}{3}}\Vert\big( D^{N-1}_{y}w_{k}(s) , D^{N-1}_{y}j_{k}(s) \big)\Vert_{L_{y}^{2}}^{2}e^{-\frac{1}{8(1+M^{'}_{k})}\nu^{\frac{1}{3}}\vert k \vert^{\frac{2}{3}}(t-s)}ds\\
  & \qquad +  \int_{0}^{t} C \nu^{-\frac{1}{3}-8}\vert k \vert^{\frac{4}{3}}\Vert D^{N}_{y}\theta_{k}(s) \Vert_{L_{y}^{2}}^{2}e^{-\frac{1}{8(1+M^{'}_{k})}\nu^{\frac{1}{3}}\vert k \vert^{\frac{2}{3}}(t-s)} ds .
\end{align*}
Next, we prove the  inequality  (\ref{the_linear_stability_result3}) by induction. For $N = 1$,
\begin{align*}
  &  \Vert D_{y}w_{k}(t) \Vert_{L_{y}^{2}}^{2} + \Vert D_{y}j_{k}(t) \Vert_{L_{y}^{2}}^{2} \\
  & \leq C_{1}\nu^{-4}
  \Vert\big( D_{y}w_{k}(0) , D_{y}j_{k}(0) \big)\Vert_{L_{y}^{2}}^{2}e^{-c_{1}\nu^{\frac{13}{3}}\vert k \vert^{\frac{2}{3}}t} \\
   & \qquad + C_{1}\nu^{-12}(\nu^{-1}\vert k \vert)^{\frac{2}{3}}\left(\nu^{-4} \Vert \big(w_{k}(0) , j_{k}(0)\big) \Vert_{L_{y}^{2}}^{2} + \nu^{-12}(\nu^{-1} \vert k \vert)^{\frac{2}{3}}\Vert \theta_{k}(0) \Vert_{L_{y}^{2}}^{2} \right)e^{-c_{1}\nu^{\frac{13}{3}}\vert k \vert^{\frac{2}{3}}t} \\
   & \qquad + C_{1}\nu^{-12}(\nu^{-1}\vert k \vert)^{\frac{2}{3}}\left( \Vert D_{y}\theta_{k}(0) \Vert_{L_{y}^{2}}^{2} + (\eta^{-1}\vert k \vert)^{\frac{2}{3}}\Vert \theta_{k}(0) \Vert_{L_{y}^{2}}^{2} \right)e^{-c_{1}\nu^{\frac{13}{3}}\vert k \vert^{\frac{2}{3}}t} 
\\
  & \leq C_{1}^{2} \left( \nu^{-4}\Vert\big( D_{y}w_{k}(0), D_{y}j_{k}(0)\big) \Vert_{L_{y}^{2}}^{2} + \nu^{-12}(\nu^{-1}\vert k \vert)^{\frac{2}{3}}\Vert D_{y}\theta_{k}(0) \Vert_{L_{y}^{2}}^{2}\right)e^{-2c_{1}\nu^{\frac{13}{3}}\vert k \vert^{\frac{2}{3}}t} \\
  & \qquad + C_{1}^{2}\nu^{ - 12}(\nu^{-1}\vert k \vert)^{\frac{2}{3}}\left( \nu^{-4}\Vert \big(w_{k}(0) , j_{k}(0)\big) \Vert_{L_{y}^{2}}^{2} + \nu^{-12}(\nu^{-1}\vert k \vert)^{\frac{2}{3}}\Vert \theta_{k}(0) \Vert_{L_{y}^{2}}^{2}  \right)e^{-2c_{1}\nu^{\frac{13}{3}}\vert k \vert^{\frac{2}{3}}t},
\end{align*}
where $C_{1}>0$ and $c_{1}>0$. 
That is, 
\begin{align*}
   &\Vert D_{y}w_{k}(t) \Vert_{L_{y}^{2}} + \Vert D_{y}j_{k}(t) \Vert_{L_{y}^{2}}\\
    &\leq C_{1} \left( \nu^{-2} \Vert \big(D_{y}w_{k}(0),D_{y}j_{k}(0)\big) \Vert_{L_{y}^{2}} + \nu^{-6}(\nu^{-1}\vert k \vert)^{\frac{1}{3}}\Vert D_{y}\theta_{k}(0) \Vert_{L_{y}^{2}}\right)e^{-c_{1}\nu^{\frac{13}{3}}\vert k \vert^{\frac{2}{3}}t} \\
   &\qquad + C_{1}\nu^{-6}(\nu^{-1}\vert k \vert)^{\frac{1}{3}}\left( \nu^{-2}\Vert\big( w_{k}(0),j_{k}(0) \big)\Vert_{L_{y}^{2}} + \nu^{-6}(\nu^{-1}\vert k \vert)^{\frac{1}{3}}\Vert \theta_{k}(0) \Vert_{L_{y}^{2}}  \right)e^{-c_{1}\nu^{\frac{13}{3}}\vert k \vert^{\frac{2}{3}}t}.
\end{align*}
Now assume that for $n \leq N$, there exist two positive constants $C_{n}$ and $c_{n} $ such that
\begin{equation*}
  \begin{split}
    & \Vert D^{n}_{y}w_{k}(t) \Vert_{L_{y}^{2}}^{2} + \Vert D^{n}_{y}j_{k}(t) \Vert_{L_{y}^{2}}^{2}\\
    & \leq C_{n}^{2} 
    \left( \nu^{-4}\Vert\big( D^{n}_{y}w_{k}(0) , D^{n}_{y}j_{k}(0)\big) \Vert_{L_{y}^{2}}^{2} + \nu^{-12}(\nu^{-1}\vert k \vert)^{\frac{2}{3}}\Vert D^{n}_{y}\theta_{k}(0) \Vert_{L_{y}^{2}}^{2}\right)e^{-2c_{n}\nu^{\frac{13}{3}}\vert k \vert^{\frac{2}{3}}t} \\
    & \qquad + C_{n}^{2}\nu^{-12n}(\nu^{-1}\vert k \vert)^{\frac{2n}{3}}\left( \nu^{-4}\Vert \big(w_{k}(0) , j_{k}(0) \big)\Vert_{L_{y}^{2}}^{2} + \nu^{-12}(\nu^{-1}\vert k \vert)^{\frac{2}{3}}\Vert \theta_{k}(0) \Vert_{L_{y}^{2}}^{2}  \right)e^{-2c_{n}\nu^{\frac{13}{3}}\vert k \vert^{\frac{2}{3}}t}.
  \end{split}
\end{equation*}
Then for $n = N+1$, we have
\begin{align*}
  &  \Vert D^{N+1}_{y}w_{k}(t) \Vert_{L_{y}^{2}}^{2} + \Vert D^{N+1}_{y}j_{k}(t) \Vert_{L_{y}^{2}}^{2} \\
  & \leq  C\nu^{-4}\left(\Vert D^{N+1}_{y}w_{k}(0) \Vert_{L_{y}^{2}}^{2} + \Vert D^{N+1}_{y}j_{k}(0) \Vert_{L_{y}^{2}}^{2}\right)e^{-\frac{1}{8(1+M^{'}_{k})}\nu^{\frac{1}{3}}\vert k \vert^{\frac{2}{3}}t} \\
  & \qquad +  \int_{0}^{t} C (N+1)^{2} \nu^{-\frac{1}{3}-8}\vert k \vert^{\frac{4}{3}}
\Vert \big(D^{N}_{y}w_{k}(s) , D^{N}_{y}j_{k}(s) \big)\Vert_{L_{y}^{2}}^{2}e^{-\frac{1}{8(1+M^{'}_{k})}\nu^{\frac{1}{3}}\vert k \vert^{\frac{2}{3}}(t-s)}ds\\
  & \qquad +  \int_{0}^{t} C \nu^{-\frac{1}{3}-8}\vert k \vert^{\frac{4}{3}}\Vert D^{N+1}_{y}\theta_{k}(s) \Vert_{L_{y}^{2}}^{2}e^{-\frac{1}{8(1+M^{'}_{k})}\nu^{\frac{1}{3}}\vert k \vert^{\frac{2}{3}}(t-s)} ds \\
  & \leq C\nu^{-4}\left(\Vert D^{N+1}_{y}w_{k}(0) \Vert_{L_{y}^{2}}^{2} + \Vert D^{N+1}_{y}j_{k}(0) \Vert_{L_{y}^{2}}^{2}\right)e^{-c\nu^{\frac{13}{3}}\vert k \vert^{\frac{2}{3}}t}+ I_{1} + I_{2}.
\end{align*}
From the induction assumption, we have
\begin{align*}
  I_{1}& = \int_{0}^{t} C (N+1)^{2} \nu^{-\frac{1}{3}-8}\vert k \vert^{\frac{4}{3}}
\Vert \big(D^{N}_{y}w_{k}(s) , D^{N}_{y}j_{k}(s) \big)\Vert_{L_{y}^{2}}^{2}e^{-\frac{1}{8(1+M^{'}_{k})}\nu^{\frac{1}{3}}\vert k \vert^{\frac{2}{3}}(t-s)}ds\\
& \leq C^{2}_{N+1}e^{-2\beta_{N}\nu^{\frac{13}{3}}\vert k \vert^{\frac{2}{3}}t} (\nu^{-1}\vert k \vert)^{\frac{2}{3}}\nu^{-12}
  \\
 & \qquad \times \left(\nu^{-4} \Vert\big( D^{N}_{y}w_{k}(0) , D^{N}_{y}j_{k}(0)\big) \Vert_{L_{y}^{2}}^{2} + \nu^{-12}(\nu^{-1}\vert k \vert)^{\frac{2}{3}}\Vert D^{N}_{y}\theta_{k}(0) \Vert_{L_{y}^{2}}^{2} \right) \\
  & \quad + C^{2}_{N+1}e^{-2\beta_{N}\nu^{\frac{13}{3}}\vert k \vert^{\frac{2}{3}}t}\nu^{-12(N+1)}(\nu^{-1}\vert k \vert)^{\frac{2(N+1)}{3}}\\
  &\qquad\times\left( \nu^{-4}\Vert (w_{k}(0),j_{k}(0)) \Vert_{L_{y}^{2}}^{2} + \nu^{-12}(\nu^{-1}\vert k \vert)^{\frac{2}{3}}\Vert \theta_{k}(0) \Vert_{L_{y}^{2}}^{2} \right),
\end{align*}
with $\beta_{N} > 0$. Note that
\begin{align*}
  &\nu^{-12}(\nu^{-1}\vert k \vert)^{\frac{2}{3}}\Vert D^{N}_{y}w_{k}(0) \Vert_{L_{y}^{2}}^{2} \leq  \nu^{-12}(\nu^{-1}\vert k \vert)^{\frac{2}{3}}\Vert w_{k}(0) \Vert_{L_{y}^{2}}^{\frac{2}{N+1}}\Vert D^{N+1}_{y}w_{k}(0) \Vert_{L_{y}^{2}}^{\frac{2N}{N+1}} \\
  &\leq \frac{1}{N+1}\nu^{-12(N+1)}(\nu^{-1}\vert k \vert)^{\frac{2(N+1)}{3}}\Vert w_{k}(0) \Vert_{L_{y}^{2}}^{2} + \frac{N}{N+1}\Vert D^{N+1}_{y}w_{k}(0) \Vert_{L_{y}^{2}}^{2}, 
  \end{align*}
  \begin{align*}
  &\nu^{-12}(\nu^{-1}\vert k \vert)^{\frac{2}{3}}\Vert D^{N}_{y}j_{k}(0) \Vert_{L_{y}^{2}}^{2} \leq \nu^{-12}(\nu^{-1}\vert k \vert)^{\frac{2}{3}}\Vert j_{k}(0) \Vert_{L_{y}^{2}}^{\frac{2}{N+1}}\Vert D^{N+1}_{y}j_{k}(0) \Vert_{L_{y}^{2}}^{\frac{2N}{N+1}} \\
 & \leq  \frac{1}{N+1}\nu^{-12(N+1)}(\nu^{-1}\vert k \vert)^{\frac{2(N+1)}{3}}\Vert j_{k}(0) \Vert_{L_{y}^{2}}^{2} + \frac{N}{N+1}\Vert D^{N+1}_{y}j_{k}(0) \Vert_{L_{y}^{2}}^{2}, 
 \end{align*}
  \begin{align*}
  &\nu^{-12}(\nu^{-1}\vert k \vert)^{\frac{2}{3}}\Vert D^{N}_{y}\theta_{k}(0) \Vert_{L_{y}^{2}}^{2} \leq  \nu^{-12}(\nu^{-1}\vert k \vert)^{\frac{2}{3}}\Vert \theta_{k}(0) \Vert_{L_{y}^{2}}^{\frac{2}{N+1}}\Vert D^{N+1}_{y}\theta_{k}(0) \Vert_{L_{y}^{2}}^{\frac{2N}{N+1}} \\
   &\leq  \frac{1}{N+1}\nu^{-12(N+1)}(\nu^{-1}\vert k \vert)^{\frac{2(N+1)}{3}}\Vert \theta_{k}(0) \Vert_{L_{y}^{2}}^{2} + \frac{N}{N+1}\Vert D^{N+1}_{y}\theta_{k}(0) \Vert_{L_{y}^{2}}^{2}.
\end{align*}
For $I_{2}$, we have
\begin{align*}
  I_{2}&=\int_{0}^{t} C \nu^{-\frac{1}{3}-8}\vert k \vert^{\frac{4}{3}}\Vert D^{N+1}_{y}\theta_{k}(s) \Vert_{L_{y}^{2}}^{2}e^{-\frac{1}{8(1+M^{'}_{k})}\nu^{\frac{1}{3}}\vert k \vert^{\frac{2}{3}}(t-s)} ds \\
   & \leq C^{2}_{N+1} e^{-2\gamma_{N}\nu^{\frac{13}{3}}\vert k \vert^{\frac{2}{3}}t}\nu^{-12} (\nu^{-1}\vert k \vert)^{\frac{2}{3}} \left( \Vert D_{y}^{N+1}\theta_{k}(0) \Vert_{L_{y}^{2}}^{2} + (\eta^{-1}\vert k \vert)^{\frac{2(N+1)}{3}}\Vert \theta_{k}(0) \Vert_{L_{y}^{2}}^{2} \right), \\
  & \leq C^{2}_{N+1} e^{-2\gamma_{N}\nu^{\frac{13}{3}}\vert k \vert^{\frac{2}{3}}t}\nu^{-12} (\nu^{-1}\vert k \vert)^{\frac{2}{3}} \left( \Vert D_{y}^{N+1}\theta_{k}(0) \Vert_{L_{y}^{2}}^{2} + (\nu^{-1}\vert k \vert)^{\frac{2(N+1)}{3}}\Vert \theta_{k}(0) \Vert_{L_{y}^{2}}^{2} \right),
\end{align*}
where the positive constants $\beta_{N}$ and $ \gamma_{N}  $ are dependent of  $N $.
Thus we have
\begin{equation*}
  \begin{split}
  & \Vert\big  (D^{N+1}_{y}w_{k}(t),  D^{N+1}_{y}j_{k}(t)\big) \Vert_{L_{y}^{2}} \\
   &\leq  C_{N+1}
e^{-c_{N+1}\nu^{\frac{13}{3}}\vert k \vert^{\frac{2}{3}}t}
\Big( \nu^{-2}\Vert\left( D^{N+1}_{y}w_{k}(0),D^{N+1}_{y}j_{k}(0) \right)\Vert_{L_{y}^{2}}  + \nu^{-6}(\nu^{-1}\vert k \vert)^{\frac{1}{3}}\Vert D^{N+1}_{y}\theta_{k}(0) \Vert_{L_{y}^{2}}  \\
&\qquad
 + \nu^{-6(N+1)}(\nu^{-1}\vert k \vert)^{\frac{N+1}{3}}
 \big( \nu^{-2}\Vert\big( w_{k}(0), j_{k}(0)\big) \Vert_{L_{y}^{2}} 
 + \nu^{-6}(\nu^{-1}\vert k \vert)^{\frac{1}{3}}\Vert \theta_{k}(0) \Vert_{L_{y}^{2}}\big) \Big),
  \end{split}
\end{equation*}
with $c_{N+1} = \min\{ \beta_{N}, \gamma_{N} \}$.

This completes the proof.
\end{proof}

\subsection{Proof of Theorem \ref{the_physics_linear_stability_result}}

In this subsection, we prove Theorem \ref{the_physics_linear_stability_result}, which is a consequence of Theorem \ref{the_linear_stability_result}.
\begin{proof}
For any $b\in \mathbb{R}$, we apply $\Lambda^{b}_{t}$ to the equations in (\ref{linear_system__projecting_y}) to obtain
\begin{equation*}
\begin{cases}
  \partial_{t}\Lambda^{b}_{t}w_{k} + iky\Lambda^{b}_{t}w_{k} +   \nu (D^{2}_{y}+ k^{2}) \Lambda^{b}_{t}w_{k} = ik \Lambda^{b}_{t}\theta_{k} + ik\Lambda^{b}_{t}j_{k}, \\ 
  \partial_{t}\Lambda^{b}_{t}j_{k} + iky\Lambda^{b}_{t}j_{k} + \nu (D^{2}_{y}+ k^{2}) \Lambda^{b}_{t}j_{k} - 2ik\Lambda^{b}_{t}b^{1}_{k}  = ik\Lambda^{b}_{t}w_{k}.
\end{cases}
\end{equation*}
We  multiply the above equations by  $(1+M^{'}_{k})\Lambda^{b}_{t}w_{k}$ and $(1+M^{'}_{k})\Lambda^{b}_{t}j_{k}$, respectively, and integrate over $\mathbb{R}$. Then using the properties of $M^{'}_{k}$, we can show that
\begin{align*}
   & \frac{d}{dt}\langle (1+M^{'}_{k})\Lambda^{b}_{t}w_{k},\Lambda^{b}_{t}w_{k} \rangle_{L_{y}^{2}} + \nu \Vert D_{y}\Lambda^{b}_{t}w_{k} \Vert^{2}_{L_{y}^{2}} +  \nu k^{2} \Vert\Lambda^{b}_{t} w_{k} \Vert^{2}_{L_{y}^{2}} + \frac{1}{8}\nu^{\frac{1}{3}}\vert k \vert^{\frac{2}{3}}\Vert\Lambda^{b}_{t} w_{k} \Vert^{2}_{L_{y}^{2}} \\
  & \qquad + \frac{d}{dt}\langle (1+M^{'}_{k})\Lambda^{b}_{t}j_{k},\Lambda^{b}_{t}j_{k} \rangle_{L_{y}^{2}} + \nu \Vert D_{y}\Lambda^{b}_{t}j_{k} \Vert^{2}_{L_{y}^{2}} +  \nu k^{2} \Vert \Lambda^{b}_{t}j_{k} \Vert^{2}_{L_{y}^{2}} + \frac{1}{8}\nu^{\frac{1}{3}}\vert k \vert^{\frac{2}{3}}\Vert \Lambda^{b}_{t}j_{k} \Vert^{2}_{L_{y}^{2}} \\
 & \leq  C\nu^{-\frac{1}{3}-8} \vert k \vert^{\frac{4}{3}}\Vert \Lambda^{b}_{t}\theta_{k} \Vert^{2}_{L_{y}^{2}}.
\end{align*}
Summing over $k$ and integrating in t yields
\begin{align*}
 & \Vert \Lambda_{t}^{b}w \Vert_{L_{t}^{\infty}(L^{2})}^{2}  - \sum_{k\in\mathbb{Z}}\langle (1+M^{'}_{k})\Lambda^{b}_{t=0}w_{k},\Lambda^{b}_{t=0}w_{k} \rangle_{L_{y}^{2}} + \nu \int_{0}^{\infty}\sum_{k\in\mathbb{Z}}\Vert D_{y}\Lambda^{b}_{t}w_{k}(t) \Vert^{2}_{L_{y}^{2}}dt \\
 & \quad +  \nu \int_{t=0}^{\infty}\sum_{k\in\mathbb{Z}} k^{2} \Vert\Lambda^{b}_{t} w_{k} \Vert^{2}_{L_{y}^{2}}dt + \frac{1}{8}\nu^{\frac{1}{3}}\int_{t=0}^{\infty}\sum_{k\in\mathbb{Z}}\vert k \vert^{\frac{2}{3}}\Vert\Lambda^{b}_{t} w_{k} \Vert^{2}_{L_{y}^{2}}dt \\
  & \quad + \Vert \Lambda_{t}^{b}j \Vert_{L_{t}^{\infty}(L^{2})}^{2} -\sum_{k\in\mathbb{Z}} \langle (1+M^{'}_{k})\Lambda^{b}_{t=0}j_{k},\Lambda^{b}_{t=0}j_{k} \rangle_{L_{y}^{2}} + \nu \int_{t=0}^{\infty}\sum_{k\in\mathbb{Z}}\Vert D_{y}\Lambda^{b}_{t}j_{k} \Vert^{2}_{L_{y}^{2}}dt \\
  &\quad  +  \nu \int_{t=0}^{\infty}\sum_{k\in\mathbb{Z}} k^{2} \Vert \Lambda^{b}_{t}j_{k} \Vert^{2}_{L_{y}^{2}}dt + \frac{1}{8}\nu^{\frac{1}{3}}\int_{t=0}^{\infty}\sum_{k\in\mathbb{Z}}\vert k \vert^{\frac{2}{3}}\Vert \Lambda^{b}_{t}j_{k} \Vert^{2}_{L_{y}^{2}}dt \\
 & \leq  C\nu^{-\frac{1}{3}-8} \int_{t=0}^{\infty}\sum_{k\in\mathbb{Z}} \vert k \vert^{\frac{4}{3}}\Vert \Lambda^{b}_{t}\theta_{k} \Vert^{2}_{L_{y}^{2}}dt.
\end{align*}
Thus we have
\begin{equation}\label{physical_inequation_in_1}
\begin{split}
  & \Vert \Lambda_{t}^{b}w \Vert_{L_{t}^{\infty}(L^{2})} + \Vert \Lambda_{t}^{b}j \Vert_{L_{t}^{\infty}(L^{2})} + \nu^{\frac{1}{2}} (\Vert \nabla \Lambda^{b}_{t}w \Vert_{L_{t}^{2}(L^{2})}+ \Vert \nabla \Lambda^{b}_{t}j \Vert_{L_{t}^{2}(L^{2})}) \\
 & \quad + \nu^{\frac{1}{6}}(\Vert \vert D_{x} \vert^{\frac{1}{3}} \Lambda^{b}_{t}w \Vert_{L_{t}^{2}(L^{2})} +  \Vert \vert D_{x} \vert^{\frac{1}{3}} \Lambda^{b}_{t}j \Vert_{L_{t}^{2}(L^{2})}) \\
 & \leq C(\nu^{-\frac{1}{6}-4} \Vert \vert D_{x} \vert^{\frac{2}{3}} \Lambda^{b}_{t}\theta \Vert_{L_{t}^{2}(L^{2})} + \nu^{-2}\Vert w(0) \Vert_{H^{b}} + \nu^{-2}\Vert j(0) \Vert_{H^{b}}).
\end{split}
\end{equation}
Similarly for $\vert D_{x} \vert^{\frac{1}{3}}\Lambda_{t}^{b}\theta$, we have
\begin{equation}\label{physical_inequation_in_2}
\Vert \Lambda_{t}^{b}\vert D_{x} \vert^{\frac{1}{3}}\theta \Vert_{L_{t}^{\infty}(L^{2})} + \eta^{\frac{1}{2}} \Vert \nabla \vert D_{x} \vert^{\frac{1}{3}}\Lambda^{b}_{t}\theta \Vert_{L_{t}^{2}(L^{2})}  + \eta^{\frac{1}{6}}\Vert \vert \vert D_{x}\vert^{\frac{2}{3}} \Lambda^{b}_{t}\theta \Vert_{L_{t}^{2}(L^{2})} \leq C\Vert \vert D_{x} \vert^{\frac{1}{3}}\theta(0) \Vert_{H^{b}}.
\end{equation}
According to (\ref{physical_inequation_in_1}) and (\ref{physical_inequation_in_2}), we get
\begin{align*}
    & \Vert \Lambda_{t}^{b}w \Vert_{L_{t}^{\infty}(L^{2})} + \nu^{\frac{1}{2}}\Vert \nabla \Lambda_{t}^{b}w \Vert_{L_{t}^{2}(L^{2})} + \nu^{\frac{1}{6}}\Vert \vert D_{x}\vert^{\frac{1}{3}} \Lambda_{t}^{b}w \Vert_{L_{t}^{2}(L^{2})} \\
    &\quad + \Vert \Lambda_{t}^{b}j \Vert_{L_{t}^{\infty}(L^{2})} + \nu^{\frac{1}{2}}\Vert \nabla\Lambda_{t}^{b}j \Vert_{L_{t}^{2}(L^{2})} + \nu^{\frac{1}{6}}\Vert \vert D_{x}\vert^{\frac{1}{3}} \Lambda_{t}^{b}j \Vert_{L_{t}^{2}(L^{2})} \\
    &\quad + \nu^{-4}(\nu\eta)^{-\frac{1}{6}}\left( \Vert \vert D_{x} \vert^{\frac{1}{3}}\Lambda_{t}^{b}\theta \Vert_{L_{t}^{\infty}(L^{2})} + \eta^{\frac{1}{2}}\Vert \nabla \vert D_{x} \vert^{\frac{1}{3}}\Lambda_{t}^{b}\theta \Vert_{L_{t}^{2}(L^{2})} + \eta^{\frac{1}{6}}\Vert\vert D_{x} \vert^{\frac{2}{3}} \Lambda_{t}^{b}\theta \Vert_{L_{t}^{2}(L^{2})} \right) \\
    & \leq C\left( \nu^{-2}\Vert w(0) \Vert_{H^{b}} + \nu^{-2}\Vert j(0) \Vert_{H^{b}} + \nu^{-4}(\nu\eta)^{-\frac{1}{6}} \Vert \vert D_{x} \vert^{\frac{1}{3}}\theta(0) \Vert_{H^{b}} \right),
\end{align*}
which concludes the proof.
\end{proof}

\section{The nonlinear stability }
For the proof of the nonlinear stability, we recall  the time-dependent elliptic operator which is defined in introduction. For $t \geq 0$,
\begin{equation*}
  \Lambda^{2}_{t} = 1 - \partial^{2}_{x} - (\partial_{y} + t\partial_{x})^{2},
\end{equation*}
or, in terms of its symbol, $\Lambda^{2}_{t}(k,\xi) = 1 + k^{2} + (\xi + tk)^{2}$.  And for any $b \in \mathbb{R}$, $\Lambda^{b}_{t}(k,\xi) = (1 + k^{2} + (\xi + tk)^{2})^{\frac{b}{2}}$ and the operator $\Lambda_{t}^{b}$ satisfies the following basic properties, which can be found in \cite{deng2020stability}.
\begin{Lemma}\label{lemma_Lambda_t_b}
  The operator $\Lambda_{t}^{b}$ satisfies the following properties
  \begin{enumerate}
  	\item For any $b \in \mathbb{R}$, $\Lambda_{t}^{b}$ commutes with $\partial_{t} + y\partial_{x}$, namely
  \[
    \Lambda_{t}^{b}(\partial_{t} + y\partial_{x}) = (\partial_{t} + y\partial_{x})\Lambda_{t}^{b}.
  \] 
  	\item For any $b>0$,
  	\[
  	  \Vert \Lambda_{t}^{b}(fg) \Vert_{L^{2}} \leq \Vert f \Vert_{L^{\infty}}\Vert \Lambda_{t}^{b}g \Vert_{L^{2}} + \Vert g \Vert_{L^{\infty}}\Vert \Lambda_{t}^{b}f \Vert_{L^{2}}.
  	\]
  	Moreover, for $b>1$, we have
  	\[
  	  \Vert f(t) \Vert_{L^{\infty}(\Omega)} \leq C \Vert \hat{f}(t) \Vert_{L^{1}} \leq C \Vert \Lambda_{t}^{b}f(t) \Vert_{L^{2}(\Omega)}.
  	\]
  	And consequently,
  	\[
  	  \Vert \Lambda_{t}^{b}(fg) \Vert_{L^{2}} \leq C \Vert \Lambda_{t}^{b}f \Vert_{L^{2}}\Vert \Lambda_{t}^{b}g \Vert_{L^{2}}.
  	\]
  \end{enumerate}
\end{Lemma}
 The framework is the bootstrap argument \cite{tao2006nonlinear}, which is stated as follows.
\begin{Lemma}[Abstract bootstrap principle] Let $I$ be a time interval, and for each $t \in I$ suppose we have two statements,  that is, ``hypothesis" $\mathbf{H}(t)$ and ``conclusion" $\mathbf{C}(t)$. Suppose we can verify the following four assertions:
\begin{enumerate}
\item (Hypothesis implies conclusion) If $\mathbf{H}(t)$ is true for some time $t\in I$, then $\mathbf{C}(t)$ is also true for that time $t$.
\item (Conclusion is stronger than hypothesis) If $\mathbf{C}(t)$ is true for some $t \in I$, then $\mathbf{H}(t^{'})$ is true for all $t^{'} \in I$ in a neighbourhood of $t$.
\item (Conclusion is closed) If $t_{1}, t_{2},...$ is a sequence of times in $I$ which converges to another time $t \in I$, and $\mathbf{C}(t_{n})$ is true for all $t_{n}$, then $\mathbf{C}(t)$ is true.
\item (Base case) $\mathbf{H}(t)$ is true for at least one time $t \in I$.
\end{enumerate}
Then $\mathbf{C}(t)$ is true for all $t \in I$.
\end{Lemma}
In order to use abstract bootstrap principle to prove Theorem \ref{nolinear_theorem}, we define $I=[0,\infty]$, the statement of $\mathbf{H}(T)$ as an estimation less than or equal to $C\varepsilon \nu^{\alpha}$ and  the statement of $\mathbf{C}(T)$ as the same estimation less than or equal to $\frac{1}{2}C\varepsilon \nu^{\alpha}$. According to non-decreasing estimation for time $t$ it is easy to check the assertion 2. According to continuity of estimation for time $t$ we get the assertion 3. And when $\mathbf{H}(t)$ is true on $t = 0$, the remaining task is to prove $\mathbf{H}(T) \Rightarrow \mathbf{C}(T)$ based on the priori bounds. Then we prove that $\mathbf{C}(t)$ is true for all $t \in [0,\infty]$. 

Let $\nu = \eta = \mu$ in (\ref{nonlinear_stability_system_rewrite}). Invoking the properties of $\Lambda_{t}^{b}$ in Lemma \ref{lemma_Lambda_t_b}, we have
\begin{equation}\label{nonlinear_stability_system_lambdatb}
\begin{cases}
  \partial_{t}\Lambda^{b}_{t}w + y\partial_{x}\Lambda^{b}_{t}w - \nu \Delta \Lambda^{b}_{t}w + \Lambda^{b}_{t}\left((\mathbf{u} \cdot \nabla) w\right) - \Lambda^{b}_{t}\left((\mathbf{b} \cdot \nabla)j\right) = \partial_{x}\Lambda^{b}_{t}j  + \partial_{x} \Lambda^{b}_{t}\theta, \\
  \partial_{t}\Lambda^{b}_{t}j + y\partial_{x}\Lambda^{b}_{t}j - \nu \Delta \Lambda^{b}_{t}j + \Lambda^{b}_{t}\left((\mathbf{u} \cdot \nabla) j\right)       \\
\quad \qquad \qquad \qquad- \Lambda^{b}_{t}\left((\mathbf{b} \cdot \nabla) w \right)
  -  2\partial_{x}\Lambda^{b}_{t}b^{1} - \Lambda^{b}_{t}Q(\nabla\mathbf{u}, \nabla \mathbf{b}) = \partial_{x} \Lambda^{b}_{t}w, \\
  \partial_{t}\Lambda^{b}_{t}\theta + y\partial_{x} \Lambda^{b}_{t}\theta - \nu\Delta\Lambda^{b}_{t}\theta + \Lambda^{b}_{t}\left((\mathbf{u}\cdot \nabla)\theta\right) = 0, \\
  \mathbf{u} = -\nabla^{\bot}(-\Delta)^{-1}w, \\
  \mathbf{b} = -\nabla^{\bot}(-\Delta)^{-1}j,
\end{cases}
\end{equation}
where $Q(\nabla\mathbf{u}, \nabla \mathbf{b}) = 2\partial_{x}b^{1}(\partial_{x}u^{2} + \partial_{y}u^{1}) - 2\partial_{x}u^{1}(\partial_{x}b^{2} + \partial_{y}b^{1}).$

Now we are ready to prove Theorem \ref{nolinear_theorem}.

{\bf Proof of Theorem \ref{nolinear_theorem}}\,\,\,
 Let $\varphi$ be a real-valued, non-decreasing function, and $\varphi \in C^{\infty}(\mathbb{R})$ satisfies $0 \leq \varphi(x) \leq 1$, $0 \leq \varphi^{'}(x) \leq \frac{1}{4}$ for all  $x \in \mathbb{R}$ and $\varphi^{'}(x) = \frac{1}{4}$ for  $x \in [-1,1]$.
 For $k \neq 0$, we choose the function $\phi_{k}$ as follows, 
\begin{equation*}
\phi_{k}(\xi) = 
\begin{cases}
  \frac{6(k^{2} + \xi_{0}^{2})^{2}}{k^{4}} - (2 + \pi) ,\qquad & \xi > 0, \\
  \frac{6(k^{2} + \xi_{0}^{2})^{2}}{(k^{2} + \xi^{2})^{2}} - (2 + \pi), \qquad & \xi \in [-\xi_{0},0], \\
  (4 - \pi) e^{\frac{24\xi_{0}(\xi+\xi_{0})}{(4-\pi)(k^{2}+\xi_{0}^{2})}}, \qquad & \xi \in (-\infty, -\xi_{0}),
\end{cases}
\end{equation*}
where $\xi_{0}$ is a real positive solution of the equation $\nu \xi_{0} (k^{2} + \xi_{0}^{2}) = 96 \vert k \vert$. 
 $\phi_{k} (\xi)\in C^{1}(\mathbb{R})$ for  $k \neq 0 $,
 and $0 < \phi_{k}(\xi) \leq C\nu^{-4}$, $0 \leq \phi_{k}^{'}(\xi) \leq C\frac{\nu^{-3}}{\vert k \vert}$, $\xi\in \mathbb{R}$
 where $C$ is a positive constant.
The Fourier multiplier operator $\mathcal{M}$ employed here is defined as $\mathcal{M} = \mathcal{M}_{1} + \mathcal{M}_{2} + \mathcal{M}_{3} + 1$ with $\mathcal{M}_{1}$, $\mathcal{M}_{2}$ and $\mathcal{M}_{3}$ given by
\begin{equation*}
\mathcal{M}_{1}(k,\xi) = \varphi(\nu^{\frac{1}{3}}\vert k \vert^{-\frac{1}{3}}{\rm sgn}(k)\xi),\, k \neq 0, 
\end{equation*}
\begin{equation*}
\mathcal{M}_{2}(k,\xi) = \phi_{k}({\rm sgn}(k)\xi),\,\quad\qquad k \neq 0,
\end{equation*}
\begin{equation*}
\mathcal{M}_{3}(k,\xi)  = \frac{1}{k^{2}}\Big({\rm  arctan} \frac{\xi}{k} + \frac{\pi}{2} \Big),\, k \neq 0, 
\end{equation*}
\begin{equation*}
  \mathcal{M}_{1}(0,\xi) = \mathcal{M}_{2}(0,\xi) = \mathcal{M}_{3}(0,\xi) = 0.
\end{equation*}
Then $\mathcal{M}$ is a self-adjoint Fourier multiplier and verifies that
\[
  1 \leq \mathcal{M} \leq (C+6)\nu^{-4}.
\]
Finally we recall the projectors onto the horizontal zeroth mode and the non-zeroth modes,
\begin{equation}\label{f_neq}
\begin{split}
&  f_{0}(y) := (\mathbb{P}_{0}f)(y) = \frac{1}{2\pi}\int_{\mathbb{T}} f(x,y) dx,\\
& f_{\neq}(x, y) := (\mathbb{P}_{\neq}f)(x, y) = f(x, y) - (\mathbb{P}_{0} f)(x, y).
\end{split}
\end{equation}
Then we have
\begin{align*}
 (f_{\neq})_{k}(y) = \frac{1}{2\pi}\int_{\mathbb{T}}f(x,y)e^{-ikx} dx = f_{k}(y), ~{\rm for} k \neq 0, \qquad 
  (f_{\neq})_{0}(y)  = 0.
\end{align*}
Multiplying the equations (\ref{nonlinear_stability_system_lambdatb}) by $\mathcal{M}\Lambda_{t}^{b}w$, $\mathcal{M}\Lambda_{t}^{b}j$ and $\mathcal{M}\Lambda_{t}^{b}\theta$, respectively, and integrating over $\mathbb{T}\times \mathbb{R}$,  we can prove that
\begin{align}
\begin{split}\label{theta_m_equation}
  & \frac{d}{dt}\Vert \sqrt{\mathcal{M}}\Lambda_{t}^{b}\theta \Vert^{2}_{L^{2}} + 2\nu\Vert \nabla\sqrt{\mathcal{M}}\Lambda_{t}^{b}\theta \Vert^{2}_{L^{2}} + \langle (k\partial_{\xi}\mathcal{M})(D)\Lambda_{t}^{b}\theta, \Lambda_{t}^{b}\theta \rangle_{L^{2}}  \\&\quad= -2 \langle \Lambda_{t}^{b}(\mathbf{u}\cdot\nabla\theta), \mathcal{M}\Lambda_{t}^{b}\theta \rangle_{L^{2}},
\end{split} 
\end{align}
\begin{align}
\begin{split}\label{w_m_equation}
  & \frac{d}{dt}\Vert \sqrt{\mathcal{M}}\Lambda_{t}^{b}w \Vert^{2}_{L^{2}} + 2\nu\Vert \nabla\sqrt{\mathcal{M}}\Lambda_{t}^{b}w \Vert^{2}_{L^{2}} + \langle (k\partial_{\xi}\mathcal{M})(D)\Lambda_{t}^{b}w, \Lambda_{t}^{b}w \rangle_{L^{2}} \\
  & \quad= -2 \langle \Lambda_{t}^{b}(\mathbf{u}\cdot\nabla w), \mathcal{M}\Lambda_{t}^{b}w \rangle_{L^{2}} + 2 \langle \Lambda_{t}^{b}(\mathbf{b}\cdot\nabla j), \mathcal{M}\Lambda_{t}^{b}w \rangle_{L^{2}} + 2 \langle \partial_{x}\Lambda_{t}^{b}\theta, \mathcal{M}\Lambda_{t}^{b}w \rangle_{L^{2}} \\
  &\qquad+ 2 \langle \partial_{x}\Lambda_{t}^{b}j, \mathcal{M}\Lambda_{t}^{b}w \rangle_{L^{2}},
\end{split} 
\end{align}
\begin{align}
\begin{split}\label{j_m_equation}
  & \frac{d}{dt}\Vert \sqrt{\mathcal{M}}\Lambda_{t}^{b}j \Vert^{2}_{L^{2}} + 2\nu\Vert \nabla\sqrt{\mathcal{M}}\Lambda_{t}^{b}j \Vert^{2}_{L^{2}} + \langle (k\partial_{\xi}\mathcal{M})(D)\Lambda_{t}^{b}j, \Lambda_{t}^{b}j \rangle_{L^{2}} - 2 \langle 2\partial_{x}\Lambda_{t}^{b}b^{1}, \mathcal{M}\Lambda_{t}^{b}j \rangle_{L^{2}}\\
  &\quad = -2 \langle \Lambda_{t}^{b}(\mathbf{u}\cdot\nabla j), \mathcal{M}\Lambda_{t}^{b}j \rangle_{L^{2}} + 2 \langle \Lambda_{t}^{b}(\mathbf{b}\cdot\nabla w), \mathcal{M}\Lambda_{t}^{b}j \rangle_{L^{2}} + 2 \langle \partial_{x}\Lambda_{t}^{b}w, \mathcal{M}\Lambda_{t}^{b}j \rangle_{L^{2}} \\
  &\qquad+ 2 \langle \Lambda_{t}^{b}Q, \mathcal{M}\Lambda_{t}^{b}j \rangle_{L^{2}}.
\end{split}
\end{align}
Multiplying the $\theta$ equation of (\ref{nonlinear_stability_system_lambdatb})  by $\mathcal{M}\vert D_{x} \vert^{\frac{2}{3}}\Lambda_{t}^{b}\theta$ and integrating over $\mathbb{T}\times \mathbb{R}$,  we have
\begin{equation}\label{theta_dx1_3_m_equation}
\begin{split}
  & \frac{d}{dt}\Vert \sqrt{\mathcal{M}}\vert D_{x} \vert^{\frac{1}{3}}\Lambda_{t}^{b}\theta \Vert^{2}_{L^{2}} + 2\nu\Vert \nabla\sqrt{\mathcal{M}}\vert D_{x} \vert^{\frac{1}{3}}\Lambda_{t}^{b}\theta \Vert^{2}_{L^{2}} + \langle (\vert D_{x} \vert^{\frac{2}{3}}k\partial_{\xi}\mathcal{M})(D)\Lambda_{t}^{b}\theta, \Lambda_{t}^{b}\theta \rangle_{L^{2}} \\
  &\quad = -2 \langle \Lambda_{t}^{b}(\mathbf{u}\cdot\nabla\theta), \vert D_{x} \vert^{\frac{2}{3}}\mathcal{M}\Lambda_{t}^{b}\theta \rangle_{L^{2}}.
\end{split}
\end{equation}
According to the definition of $\mathcal{M}_{2}$, we can get
\begin{align*}
& k \partial_{\xi}\mathcal{M}_{2}(k,\xi) + (2 + \pi + \mathcal{M}_{2}(k,\xi))\frac{4k\xi}{k^{2} + \xi^{2}} \geq 0,\, {\rm for} \ {\rm sgn}(k)\xi \in (-\xi_{0},0), \\
&   \frac{1}{4}\nu\xi^{2} + (2 + \pi + \mathcal{M}_{2}(k,\xi))\frac{4k\xi}{k^{2} + \xi^{2}} \geq 0,\, {\rm for} \ {\rm sgn}(k)\xi \in (-\infty, -\xi_{0}].
\end{align*}
By $k \partial_{\xi}\mathcal{M}_{2}(k,\xi) \geq 0$ for all $\xi \in \mathbb{R}$,  we have
\begin{equation*}
  \frac{1}{4}\nu\xi^{2} + k \partial_{\xi}\mathcal{M}_{2}(k,\xi) + (1 + \mathcal{M}_{1}(k,\xi) + \mathcal{M}_{2}(k,\xi) + \mathcal{M}_{3}(k,\xi))\frac{4k\xi}{k^{2} + \xi^{2}} \geq 0.
\end{equation*}
Hence, we can get
\begin{align*}
& 2\nu (\xi^{2} + k^{2}) \mathcal{M}(k,\xi) + k\partial_{\xi}\mathcal{M}(k,\xi) \geq \nu(\xi^{2} + k^{2}) + \frac{1}{4}\nu^{\frac{1}{3}}\vert k \vert^{\frac{2}{3}} + \frac{1}{\xi^{2} + k^{2}}, \\
&   2\nu (\xi^{2} + k^{2}) \mathcal{M}(k,\xi) + k\partial_{\xi}\mathcal{M}(k,\xi) + \mathcal{M}(k,\xi) \frac{4k\xi}{k^{2}+\xi^{2}} \geq \nu(\xi^{2} + k^{2}) + \frac{1}{4}\nu^{\frac{1}{3}}\vert k \vert^{\frac{2}{3}} + \frac{1}{\xi^{2} + k^{2}}.
\end{align*}
Therefore,
\begin{equation*}
   2\nu \Vert \nabla \sqrt{\mathcal{M}}f \Vert_{L^{2}}^{2} + \langle (k\partial_{\xi}\mathcal{M})(D)f,f \rangle_{L^{2}} \geq \nu \Vert \nabla f \Vert^{2}_{L^{2}} + \frac{1}{4}\nu^{\frac{1}{3}}\Vert \vert D_{x} \vert^{\frac{1}{3}}f \Vert^{2}_{L^{2}} + \Vert (-\Delta)^{-\frac{1}{2}}f_{\neq} \Vert^{2}_{L^{2}}, 
  \end{equation*}
  \begin{align*}
  & 2\nu \Vert \nabla \sqrt{\mathcal{M}}f \Vert_{L^{2}}^{2} + \langle (k\partial_{\xi}\mathcal{M})(D)f,f \rangle_{L^{2}} - 2\langle 2\partial_{xy}(-\Delta)^{-1}f,\mathcal{M}(D)f \rangle_{L^{2}} \\
  & \qquad \qquad\qquad \qquad\qquad \qquad\qquad \qquad
   \geq \nu \Vert \nabla f \Vert^{2}_{L^{2}} + \frac{1}{4}\nu^{\frac{1}{3}}\Vert \vert D_{x} \vert^{\frac{1}{3}}f \Vert^{2}_{L^{2}} + \Vert (-\Delta)^{-\frac{1}{2}}f_{\neq} \Vert^{2}_{L^{2}},
\end{align*}
where $f_{\neq}$ is defined in (\ref{f_neq}). 
The formulas (\ref{theta_m_equation}), (\ref{w_m_equation}) and (\ref{j_m_equation})  become
\begin{align}
\begin{split}
  & \frac{d}{dt}\Vert \sqrt{\mathcal{M}}\Lambda_{t}^{b}\theta \Vert^{2}_{L^{2}} + \nu\Vert \nabla\Lambda_{t}^{b}\theta \Vert^{2}_{L^{2}} + \frac{1}{4}\nu^{\frac{1}{3}}\Vert \vert D_{x} \vert^{\frac{1}{3}}\Lambda_{t}^{b}\theta \Vert^{2}_{L^{2}}+ \Vert ( -\Delta )^{-\frac{1}{2}}\Lambda_{t}^{b}\theta_{\neq} \Vert^{2}_{L^{2}} \\
  & \leq  -2 \underbrace{\langle \Lambda_{t}^{b}(\mathbf{u}\cdot\nabla\theta), \mathcal{M}\Lambda_{t}^{b}\theta \rangle_{L^{2}}}_{I_{1}},
\end{split} \label{bian-1}
\end{align}
\begin{align}
\begin{split}
  & \frac{d}{dt}\Vert \sqrt{\mathcal{M}}\Lambda_{t}^{b}w \Vert^{2}_{L^{2}} + \nu\Vert \nabla\Lambda_{t}^{b}w \Vert^{2}_{L^{2}} + \frac{1}{4}\nu^{\frac{1}{3}}\Vert \vert D_{x} \vert^{\frac{1}{3}}\Lambda_{t}^{b}w \Vert^{2}_{L^{2}}+ \Vert ( -\Delta )^{-\frac{1}{2}}\Lambda_{t}^{b}w_{\neq} \Vert^{2}_{L^{2}} \\
  & \leq -2 \underbrace{\langle \Lambda_{t}^{b}(\mathbf{u}\cdot\nabla w), \mathcal{M}\Lambda_{t}^{b}w \rangle_{L^{2}}}_{I_{2}} + 2 \underbrace{\langle \Lambda_{t}^{b}(\mathbf{b}\cdot\nabla j), \mathcal{M}\Lambda_{t}^{b}w \rangle_{L^{2}}}_{I_{3}} + 2 \underbrace{\langle \partial_{x}\Lambda_{t}^{b}\theta, \mathcal{M}\Lambda_{t}^{b}w \rangle_{L^{2}}}_{I_{4}} \\
  &\quad+ 2 \underbrace{\langle \partial_{x}\Lambda_{t}^{b}j, \mathcal{M}\Lambda_{t}^{b}w \rangle_{L^{2}}}_{I_{5}},
\end{split} \label{bian-2}
\end{align}
\begin{align}
\begin{split}
  & \frac{d}{dt}\Vert \sqrt{\mathcal{M}}\Lambda_{t}^{b}j \Vert^{2}_{L^{2}} + \nu\Vert \nabla\Lambda_{t}^{b}j \Vert^{2}_{L^{2}} + \frac{1}{4}\nu^{\frac{1}{3}}\Vert \vert D_{x} \vert^{\frac{1}{3}}\Lambda_{t}^{b}j \Vert^{2}_{L^{2}}+ \Vert ( -\Delta )^{-\frac{1}{2}}\Lambda_{t}^{b}j_{\neq} \Vert^{2}_{L^{2}} \\
  & \leq  -2 \underbrace{\langle \Lambda_{t}^{b}(\mathbf{u}\cdot\nabla j), \mathcal{M}\Lambda_{t}^{b}j \rangle_{L^{2}}}_{I_{6}} + 2 \underbrace{\langle \Lambda_{t}^{b}(\mathbf{b}\cdot\nabla w), \mathcal{M}\Lambda_{t}^{b}j \rangle_{L^{2}}}_{I_{7}} + 2 \underbrace{\langle \partial_{x}\Lambda_{t}^{b}w, \mathcal{M}\Lambda_{t}^{b}j \rangle_{L^{2}}}_{I_{8}} \\
  &\quad+ 2 \underbrace{\langle \Lambda_{t}^{b}Q, \mathcal{M}\Lambda_{t}^{b}j \rangle_{L^{2}}}_{I_{9}}.
\end{split}\label{bian-3}
\end{align}
And the formula (\ref{theta_dx1_3_m_equation})  becomes
\begin{align}
 &\frac{d}{dt}\Vert \sqrt{\mathcal{M}}\vert D_{x} \vert^{\frac{1}{3}}\Lambda_{t}^{b}\theta \Vert^{2}_{L^{2}} + \nu\Vert \nabla\vert D_{x} \vert^{\frac{1}{3}}\Lambda_{t}^{b}\theta \Vert^{2}_{L^{2}} + \frac{1}{4}\nu^{\frac{1}{3}}\Vert \vert D_{x} \vert^{\frac{2}{3}}\Lambda_{t}^{b}\theta \Vert^{2}_{L^{2}}+ \Vert ( -\Delta )^{-\frac{1}{2}}\vert D_{x} \vert^{\frac{1}{3}}\Lambda_{t}^{b}\theta_{\neq} \Vert^{2}_{L^{2}} \notag\\
&\leq -2 \underbrace{\langle \Lambda_{t}^{b}(\mathbf{u}\cdot\nabla\theta), \vert D_{x} \vert^{\frac{2}{3}}\mathcal{M}\Lambda_{t}^{b}\theta \rangle_{L^{2}}}_{I_{10}}.\label{bian-4}
\end{align}
Since $\mathbf{u}$ and $\mathbf{b}$ are given by $w$ and $j$ via the Biot-Savart law,
\begin{equation*}
\mathbf{u} = - \nabla^{\perp}(-\Delta)^{-1} w = \left( \begin{aligned}
  \partial_{y}(-\Delta)^{-1}w \\
  - \partial_{x}(-\Delta)^{-1}w
\end{aligned} \right),
\end{equation*}
\begin{equation*}
\mathbf{b} = - \nabla^{\perp}(-\Delta)^{-1} j = \left( \begin{aligned}
  \partial_{y}(-\Delta)^{-1}j \\
  - \partial_{x}(-\Delta)^{-1}j
\end{aligned} \right),
\end{equation*}
we can decompose $\mathbf{u}$ into two parts according to (\ref{f_neq}),
\begin{equation*}
\begin{split}
  \mathbf{u}_{0} = \mathbb{P}_{0}\mathbf{u}  = \mathbb{P}_{0}\left( \begin{split}
  \partial_{y}(-\Delta)^{-1}w \\
  - \partial_{x}(-\Delta)^{-1}w
\end{split} \right) = \left( \begin{split}
  \frac{1}{2\pi}\int_{\mathbb{T}}\partial_{y}(-\Delta)^{-1}w(x,y) dx \\
  - \frac{1}{2\pi}\int_{\mathbb{T}}\partial_{x}(-\Delta)^{-1}w(x,y) dx
\end{split} \right).
\end{split}
\end{equation*}
Let $w(x,y) = - \Delta f(x,y) = -\partial_{x}^{2}f(x,y) - \partial_{y}^{2}f(x,y)$, then
\begin{equation*}
 w_{0} = -\frac{1}{2\pi}\int_{\mathbb{T}}(\partial_{x}^{2}f(x,y) +  \partial_{y}^{2}f(x,y))dx = -\frac{1}{2\pi}\int_{\mathbb{T}}\partial_{y}^{2}f(x,y)dx = - \partial_{y}^{2}f_{0},
\end{equation*}
\begin{equation*}
 \frac{1}{2\pi}\int_{\mathbb{T}}\partial_{y}(-\Delta)^{-1}w(x,y) dx  = \partial_{y}f_{0} = \partial_{y}(-\partial_{y}^{2})^{-1}w_{0}.
\end{equation*}
Therefore,
\[
  \mathbf{u}_{0} = \mathbb{P}_{0}\mathbf{u}  = \left( \begin{aligned}
  u_{0} \\
  0
\end{aligned} \right),\qquad {\rm with}\, \,u_{0} = \partial_{y}(-\partial_{y}^{2})^{-1}w_{0},
\]
similarly for $\mathbf{b}_{0}$,
\[
  \mathbf{b}_{0} = \mathbb{P}_{0}\mathbf{b}  = \left( \begin{aligned}
  b_{0} \\
  0
\end{aligned} \right),\qquad {\rm with}\, \,b_{0} = \partial_{y}(-\partial_{y}^{2})^{-1}j_{0}.
\]
Note that
\begin{equation*}
  \mathbf{u}_{\neq}  = \mathbb{P}_{\neq}\mathbf{u} = \mathbf{u} - \mathbf{u}_{0} = -\nabla^{\perp}(-\Delta)^{-1}w_{\neq}, 
  \end{equation*}
  \begin{equation*}
  \mathbf{j}_{\neq}  = \mathbb{P}_{\neq}\mathbf{j} = \mathbf{j} - \mathbf{j}_{0} = -\nabla^{\perp}(-\Delta)^{-1}j_{\neq},
\end{equation*}
Therefore,  we can write
\[
  I_{1} = \langle \Lambda_{t}^{b}(\mathbf{u}\cdot\nabla \theta),\mathcal{M}\Lambda_{t}^{b}\theta \rangle_{L^{2}} = I_{1_{1}} + I_{1_{2}},
\]
with
\[
  I_{1_{1}} = \langle \Lambda_{t}^{b}(\mathbf{u}_{\neq}\cdot\nabla \theta),\mathcal{M}\Lambda_{t}^{b}\theta \rangle_{L^{2}}, \qquad  I_{1_{2}} = \langle \Lambda_{t}^{b}(\mathbf{u}_{0}\cdot\nabla \theta),\mathcal{M}\Lambda_{t}^{b}\theta \rangle_{L^{2}}.
\]
Using the boundedness of $\mathcal{M}$ and Lemma \ref{lemma_Lambda_t_b}, we have for $b>1$,
\begin{align*}
  \vert I_{1_{1}} \vert & \leq C\nu^{-4} \Vert \Lambda_{t}^{b}(\mathbf{u}_{\neq}\cdot \nabla \theta) \Vert_{L^{2}}\Vert \Lambda_{t}^{b}\theta \Vert_{L^{2}} \\   
 & \leq C\nu^{-4} \Vert \Lambda_{t}^{b}\mathbf{u}_{\neq} \Vert_{L^{2}}\Vert \nabla \Lambda_{t}^{b}  \theta \Vert_{L^{2}}\Vert \Lambda_{t}^{b}\theta \Vert_{L^{2}} \\
 & \leq C\nu^{-4} \Vert \nabla^{\perp} (-\Delta)^{-1}\Lambda_{t}^{b}w_{\neq} \Vert_{L^{2}}\Vert \nabla \Lambda_{t}^{b}  \theta \Vert_{L^{2}}\Vert \Lambda_{t}^{b}\theta \Vert_{L^{2}} \\
 & \leq C\nu^{-4} \Vert (-\Delta)^{-\frac{1}{2}}\Lambda_{t}^{b}w_{\neq} \Vert_{L^{2}}\Vert \nabla \Lambda_{t}^{b}  \theta \Vert_{L^{2}}\Vert \Lambda_{t}^{b}\theta \Vert_{L^{2}}.
\end{align*}
We write $\mathcal{M}_{t}^{b} = \sqrt{\mathcal{M}}\Lambda_{t}^{b}$ or
\[
  \mathcal{M}_{t}^{b}(k,\xi) := \sqrt{\mathcal{M}(k,\xi)}(1+k^{2} + (\xi + kt)^{2})^{\frac{b}{2}}.
\]
Using the explicit expression of $\mathcal{M}_{t}^{b}$, we deduce that
\begin{equation}\label{partial_m_t_b_inequation}
  \vert \partial_{\xi}\mathcal{M}_{t}^{b}(k,\xi) \vert \leq C(\nu^{\frac{1}{3}}\vert k \vert^{-\frac{1}{3}} + \frac{\nu^{-3}}{\vert k \vert})(1+k^{2} + (\xi + kt)^{2})^{\frac{b}{2}}.
\end{equation}
Since $\theta_{0}$ is independent of $x$, we have
\begin{align*}
  I_{1_{2}} & = \langle \Lambda_{t}^{b}(\mathbf{u}_{0}\cdot\nabla \theta),\mathcal{M}\Lambda_{t}^{b}\theta \rangle_{L^{2}}  = \langle \Lambda_{t}^{b}(u_{0}\partial_{x}\theta_{\neq}),\mathcal{M}\Lambda_{t}^{b}\theta \rangle_{L^{2}}.
\end{align*}
Due to the the cancellations
\begin{equation*}
   \langle \mathcal{M}_{t}^{b}(u_{0}\partial_{x}\theta_{\neq}),\mathcal{M}_{t}^{b}\theta_{0} \rangle_{L^{2}} = \langle \mathcal{M}_{t}^{b}(u_{0}\theta_{\neq}),\partial_{x}\mathcal{M}_{t}^{b}\theta_{0} \rangle_{L^{2}} = 0,
\end{equation*}
\begin{equation*}
   \langle u_{0}\partial_{x}(\mathcal{M}_{t}^{b}\theta_{\neq}), \mathcal{M}_{t}^{b}\theta_{\neq} \rangle_{L^{2}} = 0,
\end{equation*} 
we have
\begin{align*}
  I_{1_{2}} & = \langle \mathcal{M}_{t}^{b}(u_{0}\partial_{x}\theta_{\neq}), \mathcal{M}_{t}^{b}\theta_{\neq} \rangle_{L^{2}} \\
  & = \langle \mathcal{M}_{t}^{b}(u_{0}\partial_{x}\theta_{\neq}) - u_{0}\partial_{x}(\mathcal{M}_{t}^{b}\theta_{\neq}), \mathcal{M}_{t}^{b}\theta_{\neq} \rangle_{L^{2}}.
\end{align*}
By $\widehat{fg} = \hat{f} * \hat{g}$ and Plancherel's theorem,
\begin{align*}
  I_{1_{2}} &= \sum_{k \neq 0}\int_{\mathbb{R}}\int_{\mathbb{R}} (\mathcal{M}_{t}^{b}(k,\xi) - \mathcal{M}_{t}^{b}(k,\xi - z)) \hat{u}(0,z)ik\widehat{\theta_{\neq}}(k,\xi -z)\mathcal{M}_{t}^{b}(k,\xi)\overline{\widehat{\theta_{\neq}}(k,\xi)} d\xi dz \\
 & = - \sum_{k \neq 0}\int_{\mathbb{R}}\int_{\mathbb{R}} (\mathcal{M}_{t}^{b}(k,\xi) - \mathcal{M}_{t}^{b}(k,\xi - z)) \frac{1}{z}\hat{w}(0,z)k\widehat{\theta_{\neq}}(k,\xi -z)\mathcal{M}_{t}^{b}(k,\xi)\overline{\widehat{\theta_{\neq}}(k,\xi)} d\xi dz.
\end{align*}
By Taylor's formula,
\begin{align*}
  \vert \mathcal{M}_{t}^{b}(k,\xi) - \mathcal{M}_{t}^{b}(k,\xi - z)) \vert \leq \int_{0}^{1} \vert \partial_{\xi}\mathcal{M}_{t}^{b}(k,\xi -sz) \vert \vert z \vert ds.
\end{align*}
Therefore by (\ref{partial_m_t_b_inequation}) and $\max_{s \in [0,1]}\vert \xi - sz + kt \vert \leq 2\vert z \vert + \vert \xi + kt -z \vert$, we get
\begin{align*}
  \vert \mathcal{M}_{t}^{b}(k,\xi) - \mathcal{M}_{t}^{b}(k,\xi - z)) \vert & \leq \int_{0}^{1} \vert \partial_{\xi}\mathcal{M}_{t}^{b}(k,\xi -sz) \vert \vert z \vert ds \\
  & \leq \int_{0}^{1} C(\nu^{\frac{1}{3}}\vert k \vert^{-\frac{1}{3}} + \frac{\nu^{-3}}{\vert k \vert})(1+k^{2} + (\xi - sz + kt)^{2})^{\frac{b}{2}} \vert z \vert ds \\
  & \leq \int_{0}^{1} C(\nu^{\frac{1}{3}}\vert k \vert^{-\frac{1}{3}} + \frac{\nu^{-3}}{\vert k \vert})(1+k^{2} + z^{2} + ( \xi + kt -z )^{2} )^{\frac{b}{2}} \vert z \vert ds \\
  & \leq C(\nu^{\frac{1}{3}}\vert k \vert^{-\frac{1}{3}} + \frac{\nu^{-3}}{\vert k \vert})\left((1+k^{2} + ( \xi + kt -z )^{2})^{\frac{b}{2}} + (1+z^{2})^{\frac{b}{2}} \right) \vert z \vert.
\end{align*}
So $I_{1_{2}}$ can be estimated as 
\begin{align*}
  \vert I_{1_{2}} \vert & \leq \sum_{k \neq 0} C(\nu^{\frac{1}{3}}\vert k \vert^{\frac{2}{3}} +  \nu^{-3}) \int_{\mathbb{R}}\int_{\mathbb{R}}\left( \Lambda_{t}^{b}(k,\xi - z) + \Lambda_{t}^{b}(0,z) \right)  \\
  &\qquad \qquad\qquad \qquad\qquad \qquad \quad\times\vert \hat{w}(0,z) \vert \vert \widehat{\theta_{\neq}}(k,\xi -z) \vert \mathcal{M}_{t}^{b}(k,\xi)\vert \widehat{\theta_{\neq}}(k,\xi)\vert  d\xi dz \\
  & \leq  C \nu^{-\frac{5}{3}} \Vert \widehat{w_{0}} \Vert_{L^{1}} \Vert \vert D_{x} \vert^{\frac{1}{3}} \Lambda_{t}^{b}\theta_{\neq} \Vert^{2}_{L^{2}} + C \nu^{-\frac{5}{3}} \Vert \Lambda_{t}^{b}w_{0} \Vert_{L^{2}} \Vert \widehat{\vert D_{x} \vert^{\frac{1}{3}} \theta_{\neq}} \Vert_{L^{1}} \Vert \vert D_{x} \vert^{\frac{1}{3}} \Lambda_{t}^{b}\theta_{\neq} \Vert_{L^{2}} \\
  & \qquad +  C\nu^{-5} \Vert \widehat{w_{0}} \Vert_{L^{1}} \Vert \Lambda_{t}^{b}\theta_{\neq} \Vert^{2}_{L^{2}} + C\nu^{-5} \Vert \Lambda_{t}^{b}w_{0} \Vert_{L^{2}} \Vert \widehat{\theta_{\neq}} \Vert_{L^{1}} \Vert \Lambda_{t}^{b}\theta_{\neq} \Vert_{L^{2}} \\
  & \leq C \nu^{-\frac{5}{3}} \Vert \Lambda_{t}^{b}w_{0} \Vert_{L^{2}} \Vert \vert D_{x} \vert^{\frac{1}{3}} \Lambda_{t}^{b}\theta \Vert^{2}_{L^{2}} + C\nu^{-5} \Vert \Lambda_{t}^{b}w_{0} \Vert_{L^{2}} \Vert \Lambda_{t}^{b}\theta_{\neq} \Vert^{2}_{L^{2}}.
\end{align*}
Consequently,
\begin{equation}\label{1}
\begin{split}
   \vert I_{1} \vert & \leq C\nu^{-4} \Vert (-\Delta)^{-\frac{1}{2}}\Lambda_{t}^{b}w_{\neq} \Vert_{L^{2}}\Vert \nabla \Lambda_{t}^{b}  \theta \Vert_{L^{2}}\Vert \Lambda_{t}^{b}\theta \Vert_{L^{2}} + C \nu^{-\frac{5}{3}} \Vert \Lambda_{t}^{b}w_{0} \Vert_{L^{2}} \Vert \vert D_{x} \vert^{\frac{1}{3}} \Lambda_{t}^{b}\theta \Vert^{2}_{L^{2}} \\
   & \qquad + C \nu^{-5} \Vert \Lambda_{t}^{b}w_{0} \Vert_{L^{2}} \Vert \Lambda_{t}^{b}\theta_{\neq} \Vert^{2}_{L^{2}}, \\
   & \leq C\nu^{-4} \Vert (-\Delta)^{-\frac{1}{2}}\Lambda_{t}^{b}w_{\neq} \Vert_{L^{2}}\Vert \nabla \Lambda_{t}^{b}  \theta \Vert_{L^{2}}\Vert \Lambda_{t}^{b}\theta \Vert_{L^{2}} + C \nu^{-\frac{5}{3}} \Vert \Lambda_{t}^{b}w_{0} \Vert_{L^{2}} \Vert \vert D_{x} \vert^{\frac{1}{3}} \Lambda_{t}^{b}\theta \Vert^{2}_{L^{2}} \\
   & \qquad + C\nu^{-5} \Vert \Lambda_{t}^{b}w_{0} \Vert_{L^{2}} \Vert (-\Delta)^{-\frac{1}{2}} \Lambda_{t}^{b}\theta_{\neq} \Vert_{L^{2}} \Vert \nabla \Lambda_{t}^{b}\theta \Vert_{L^{2}}.
\end{split}
\end{equation}
Similarly as $I_{1}$,  we can estimate $I_{2}$ and $I_{6}$ as follows 
\begin{align}
   \vert I_{2} \vert & \leq C\nu^{-4} \Vert (-\Delta)^{-\frac{1}{2}}\Lambda_{t}^{b}w_{\neq} \Vert_{L^{2}}\Vert \nabla \Lambda_{t}^{b}  w \Vert_{L^{2}}\Vert \Lambda_{t}^{b}w \Vert_{L^{2}}  + C \nu^{-\frac{5}{3}} \Vert \Lambda_{t}^{b}w_{0} \Vert_{L^{2}} \Vert \vert D_{x} \vert^{\frac{1}{3}} \Lambda_{t}^{b}w \Vert^{2}_{L^{2}} \notag\\
   & \qquad + C \nu^{-5}\Vert \Lambda_{t}^{b}w_{0} \Vert_{L^{2}} \Vert \Lambda_{t}^{b}w_{\neq} \Vert^{2}_{L^{2}}\notag\\
   & \leq C\nu^{-4} \Vert (-\Delta)^{-\frac{1}{2}}\Lambda_{t}^{b}w_{\neq} \Vert_{L^{2}}\Vert \nabla \Lambda_{t}^{b}  w \Vert_{L^{2}}\Vert \Lambda_{t}^{b}w \Vert_{L^{2}}  + C \nu^{-\frac{5}{3}} \Vert \Lambda_{t}^{b}w_{0} \Vert_{L^{2}} \Vert \vert D_{x} \vert^{\frac{1}{3}} \Lambda_{t}^{b}w \Vert^{2}_{L^{2}}   \label{2}\\
   & \qquad + C\nu^{-5} \Vert \Lambda_{t}^{b}w_{0} \Vert_{L^{2}} \Vert (-\Delta)^{-\frac{1}{2}} \Lambda_{t}^{b}w_{\neq} \Vert_{L^{2}} \Vert \nabla \Lambda_{t}^{b}w \Vert_{L^{2}},
\notag 
\end{align}
\begin{align}
   \vert I_{6} \vert & \leq C\nu^{-4} \Vert (-\Delta)^{-\frac{1}{2}}\Lambda_{t}^{b}w_{\neq} \Vert_{L^{2}}\Vert \nabla \Lambda_{t}^{b} j \Vert_{L^{2}}\Vert \Lambda_{t}^{b}j \Vert_{L^{2}}  + C \nu^{-\frac{5}{3}} \Vert \Lambda_{t}^{b}w_{0} \Vert_{L^{2}} \Vert \vert D_{x} \vert^{\frac{1}{3}} \Lambda_{t}^{b}j \Vert^{2}_{L^{2}} \notag \\
   &\qquad + C\nu^{-5} \Vert \Lambda_{t}^{b}w_{0} \Vert_{L^{2}} \Vert \Lambda_{t}^{b}j_{\neq} \Vert^{2}_{L^{2}} \notag \\
   & \leq C \nu^{-4} \Vert (-\Delta)^{-\frac{1}{2}}\Lambda_{t}^{b}w_{\neq} \Vert_{L^{2}}\Vert \nabla \Lambda_{t}^{b} j \Vert_{L^{2}}\Vert \Lambda_{t}^{b}j \Vert_{L^{2}} 
    + C \nu^{-\frac{5}{3}} \Vert \Lambda_{t}^{b}w_{0} \Vert_{L^{2}} \Vert \vert D_{x} \vert^{\frac{1}{3}} \Lambda_{t}^{b}j \Vert^{2}_{L^{2}}  \label{6}\\
   & \qquad + C \nu^{-5}\Vert \Lambda_{t}^{b}w_{0} \Vert_{L^{2}} \Vert (-\Delta)^{-\frac{1}{2}} \Lambda_{t}^{b}j_{\neq} \Vert_{L^{2}} \Vert \nabla \Lambda_{t}^{b}j \Vert_{L^{2}}.
\notag
\end{align}
Due to the cancellations
\begin{equation*}
   \langle \mathcal{M}_{t}^{b}(b_{0}\partial_{x}w_{\neq}),\mathcal{M}_{t}^{b}j_{0} \rangle_{L^{2}} = 0,
\end{equation*}
\begin{equation*}
   \langle \mathcal{M}_{t}^{b}(b_{0}\partial_{x}j_{\neq}),\mathcal{M}_{t}^{b}w_{0} \rangle_{L^{2}} = 0, 
\end{equation*}
\begin{equation*}
  \langle b_{0}\partial_{x}(\mathcal{M}_{t}^{b}w_{\neq}), \mathcal{M}_{t}^{b}j_{\neq} \rangle_{L^{2}} + \langle b_{0}\partial_{x}(\mathcal{M}_{t}^{b}j_{\neq}), \mathcal{M}_{t}^{b}w_{\neq} \rangle_{L^{2}} = 0,
\end{equation*}
 similarly as $I_{1}$, we can estimate $I_{3} + I_{7}$  as follows
\begin{align}
\vert I_{3} + I_{7} \vert & \leq C \nu^{-4} \Vert (-\Delta)^{-\frac{1}{2}}\Lambda_{t}^{b}j_{\neq} \Vert_{L^{2}}\Vert \nabla \Lambda_{t}^{b} j \Vert_{L^{2}}\Vert \Lambda_{t}^{b}w \Vert_{L^{2}} \notag\\&\quad+ C\nu^{-5} \Vert \Lambda_{t}^{b}j_{0} \Vert_{L^{2}} \Vert \Lambda_{t}^{b}j_{\neq} \Vert_{L^{2}} \Vert \Lambda_{t}^{b}w_{\neq} \Vert_{L^{2}} \notag\\
   &\quad + C \nu^{-\frac{5}{3}} \Vert \Lambda_{t}^{b}j_{0} \Vert_{L^{2}} \Vert \vert D_{x} \vert^{\frac{1}{3}} \Lambda_{t}^{b}j \Vert_{L^{2}} \Vert \vert D_{x} \vert^{\frac{1}{3}} \Lambda_{t}^{b}w \Vert_{L^{2}} \notag\\
   &\quad + C\nu^{-4} \Vert (-\Delta)^{-\frac{1}{2}}\Lambda_{t}^{b}j_{\neq} \Vert_{L^{2}}\Vert \nabla \Lambda_{t}^{b} w \Vert_{L^{2}}\Vert \Lambda_{t}^{b}j \Vert_{L^{2}} 
\notag\\&\quad+ C \nu^{-5} \Vert \Lambda_{t}^{b}j_{0} \Vert_{L^{2}} \Vert \Lambda_{t}^{b}w_{\neq} \Vert_{L^{2}} \Vert \Lambda_{t}^{b}j_{\neq} \Vert_{L^{2}} \notag\\
   & \quad + C \nu^{-\frac{5}{3}} \Vert \Lambda_{t}^{b}j_{0} \Vert_{L^{2}} \Vert \vert D_{x} \vert^{\frac{1}{3}} \Lambda_{t}^{b}w \Vert_{L^{2}} \Vert \vert D_{x} \vert^{\frac{1}{3}} \Lambda_{t}^{b}j \Vert_{L^{2}} \notag\\
   & \leq C \nu^{-4} \Vert (-\Delta)^{-\frac{1}{2}}\Lambda_{t}^{b}j_{\neq} \Vert_{L^{2}}\Vert \nabla \Lambda_{t}^{b}  j \Vert_{L^{2}}\Vert \Lambda_{t}^{b}w \Vert_{L^{2}} \label{3and7}\\
   & \quad + C\nu^{-5} \Vert \Lambda_{t}^{b}j_{0} \Vert_{L^{2}} \Vert (-\Delta)^{-\frac{1}{2}} \Lambda_{t}^{b}j_{\neq} \Vert_{L^{2}}^{\frac{1}{2}} \Vert \nabla \Lambda_{t}^{b}j \Vert_{L^{2}}^{\frac{1}{2}}\Vert (-\Delta)^{-\frac{1}{2}} \Lambda_{t}^{b}w_{\neq} \Vert_{L^{2}}^{\frac{1}{2}} \Vert \nabla \Lambda_{t}^{b}w \Vert_{L^{2}}^{\frac{1}{2}}\notag\\
   &\quad + C \nu^{-\frac{5}{3}} \Vert \Lambda_{t}^{b}j_{0} \Vert_{L^{2}} \Vert \vert D_{x} \vert^{\frac{1}{3}} \Lambda_{t}^{b}j \Vert_{L^{2}} \Vert \vert D_{x} \vert^{\frac{1}{3}} \Lambda_{t}^{b}w \Vert_{L^{2}} \notag\\
   &\quad + C \nu^{-4} \Vert (-\Delta)^{-\frac{1}{2}}\Lambda_{t}^{b}j_{\neq} \Vert_{L^{2}}\Vert \nabla \Lambda_{t}^{b} w \Vert_{L^{2}}\Vert \Lambda_{t}^{b}j \Vert_{L^{2}} \notag\\
   & \quad + C\nu^{-5} \Vert \Lambda_{t}^{b}j_{0} \Vert_{L^{2}} \Vert (-\Delta)^{-\frac{1}{2}} \Lambda_{t}^{b}w_{\neq} \Vert_{L^{2}}^{\frac{1}{2}} \Vert \nabla \Lambda_{t}^{b}w \Vert_{L^{2}}^{\frac{1}{2}}\Vert (-\Delta)^{-\frac{1}{2}} \Lambda_{t}^{b}j_{\neq} \Vert_{L^{2}}^{\frac{1}{2}} \Vert \nabla \Lambda_{t}^{b}j \Vert_{L^{2}}^{\frac{1}{2}}\notag\\
   &\quad + C \nu^{-\frac{5}{3}} \Vert \Lambda_{t}^{b}j_{0} \Vert_{L^{2}} \Vert \vert D_{x} \vert^{\frac{1}{3}} \Lambda_{t}^{b}w \Vert_{L^{2}} \Vert \vert D_{x} \vert^{\frac{1}{3}} \Lambda_{t}^{b}j \Vert_{L^{2}}\notag .
\end{align}
Using the self-adjointness of $\mathcal{M}$ and the skew-adjointness of $\partial_{x}$, we have
\begin{equation}\label{58}
\begin{split}
  I_{5} + I_{8} = \langle \partial_{x}\Lambda_{t}^{b}j, \mathcal{M}\Lambda_{t}^{b}w \rangle_{L^{2}} + \langle \partial_{x}\Lambda_{t}^{b}w, \mathcal{M}\Lambda_{t}^{b}j \rangle_{L^{2}} = 0.
\end{split}
\end{equation}
Using the upper bound  of $\mathcal{M}$, we have
\begin{equation}\label{4}
\begin{split}
  \vert I_{4} \vert & = \vert \langle \partial_{x}\Lambda_{t}^{b}\theta, \mathcal{M}\Lambda_{t}^{b}w \rangle_{L^{2}}\vert \leq C\nu^{-4} \Vert \vert D_{x}\vert^{\frac{2}{3}}\Lambda_{t}^{b}\theta \Vert_{L^{2}} \Vert \vert D_{x}\vert^{\frac{1}{3}}\Lambda_{t}^{b}w \Vert_{L^{2}}.
\end{split}
\end{equation}
For $I_{9}$, we have
\begin{align*}
  I_{9} & = \langle \Lambda_{t}^{b}Q, \mathcal{M}\Lambda_{t}^{b}j \rangle_{L^{2}} \\
  &  = \langle \Lambda_{t}^{b}\left(2\partial_{x}b^{1}(\partial_{x}u^{2} + \partial_{y}u^{1}) - 2\partial_{x}u^{1}(\partial_{x}b^{2} + \partial_{y}b^{1})\right), \mathcal{M}\Lambda_{t}^{b}j \rangle_{L^{2}} \\
  &  = \langle \Lambda_{t}^{b}\left(2\partial_{x}b^{1}(2\partial_{x}u^{2} - w) - 2\partial_{x}u^{1}(2\partial_{x}b^{2} - j)\right), \mathcal{M}\Lambda_{t}^{b}j \rangle_{L^{2}}.
\end{align*}
For $I_{9_{1}} = \langle \Lambda_{t}^{b}(\partial_{x}b^{1}\partial_{x}u^{2}), \mathcal{M}\Lambda_{t}^{b}j \rangle_{L^{2}}$,
\begin{align*}
  I_{9_{1}} & \leq C\nu^{-4}\Vert \Lambda_{t}^{b}\partial_{x}b^{1} \Vert_{L^{2}} \Vert \Lambda_{t}^{b} \partial_{x}u^{2} \Vert_{L^{2}} \Vert \Lambda_{t}^{b}j \Vert_{L^{2}} \\
  & \leq C\nu^{-4}\Vert \Lambda_{t}^{b}\partial_{x}\partial_{y}(-\Delta)^{-1}j_{\neq} \Vert_{L^{2}} \Vert \Lambda_{t}^{b} \partial_{x}\partial_{x}(-\Delta)^{-1}w_{\neq} \Vert_{L^{2}}\Vert \Lambda_{t}^{b}j \Vert_{L^{2}}  \\
  & \leq C\nu^{-4}\Vert \Lambda_{t}^{b}\widehat{j_{\neq}} \Vert_{L^{2}} \Vert \Lambda_{t}^{b} \widehat{w_{\neq}} \Vert_{L^{2}}\Vert \Lambda_{t}^{b}j \Vert_{L^{2}} \\
  & \leq C\nu^{-4} \Vert (-\Delta)^{-\frac{1}{2}} \Lambda_{t}^{b}j_{\neq} \Vert_{L^{2}}^{\frac{1}{2}} \Vert \nabla \Lambda_{t}^{b}j \Vert_{L^{2}}^{\frac{1}{2}} \Vert (-\Delta)^{-\frac{1}{2}} \Lambda_{t}^{b}w_{\neq} \Vert_{L^{2}}^{\frac{1}{2}} \Vert \nabla \Lambda_{t}^{b}w \Vert_{L^{2}}^{\frac{1}{2}} \Vert \Lambda_{t}^{b}j \Vert_{L^{2}}.
\end{align*}
For $I_{9_{2}} = \langle \Lambda_{t}^{b}(\partial_{x}b^{1}w), \mathcal{M}\Lambda_{t}^{b}j \rangle_{L^{2}}$,
\begin{align*}
  I_{9_{2}} & = \langle \Lambda_{t}^{b}(w_{0}\partial_{x}b^{1}), \mathcal{M}\Lambda_{t}^{b}j \rangle_{L^{2}} + \langle \Lambda_{t}^{b}(w_{\neq}\partial_{x}b^{1}), \mathcal{M}\Lambda_{t}^{b}j \rangle_{L^{2}} \\
  & = \langle \Lambda_{t}^{b}(w_{0}\partial_{x}b^{1}), \mathcal{M}\Lambda_{t}^{b}j_{\neq} \rangle_{L^{2}} + \langle \Lambda_{t}^{b}(w_{\neq}\partial_{x}b^{1}), \mathcal{M}\Lambda_{t}^{b}j \rangle_{L^{2}} \\
  & \leq C\nu^{-4}\Vert \Lambda_{t}^{b}w_{0} \Vert_{L^{2}} \Vert \Lambda_{t}^{b} \partial_{x}b^{1} \Vert_{L^{2}} \Vert \Lambda_{t}^{b}j_{\neq} \Vert_{L^{2}} + C\nu^{-4}\Vert \Lambda_{t}^{b}w_{\neq} \Vert_{L^{2}} \Vert \Lambda_{t}^{b} \partial_{x}b^{1} \Vert_{L^{2}} \Vert \Lambda_{t}^{b}j \Vert_{L^{2}} \\  
  & \leq C \nu^{-4}\Vert \Lambda_{t}^{b} w_{0} \Vert_{L^{2}}\Vert (-\Delta)^{-\frac{1}{2}} \Lambda_{t}^{b}j_{\neq} \Vert_{L^{2}} \Vert \nabla \Lambda_{t}^{b}j \Vert_{L^{2}} \\
  & \qquad +  C\nu^{-4} \Vert (-\Delta)^{-\frac{1}{2}} \Lambda_{t}^{b}j_{\neq} \Vert_{L^{2}}^{\frac{1}{2}} \Vert \nabla \Lambda_{t}^{b}j \Vert_{L^{2}}^{\frac{1}{2}} \Vert (-\Delta)^{-\frac{1}{2}} \Lambda_{t}^{b}w_{\neq} \Vert_{L^{2}}^{\frac{1}{2}} \Vert \nabla \Lambda_{t}^{b}w \Vert_{L^{2}}^{\frac{1}{2}} \Vert \Lambda_{t}^{b}j \Vert_{L^{2}}.
\end{align*}
For  $I_{9_{3}} = \langle \Lambda_{t}^{b}(\partial_{x}u^{1}\partial_{x}b^{2}), \mathcal{M}\Lambda_{t}^{b}j \rangle_{L^{2}}$,
\begin{align*}
  I_{9_{3}} & \leq C\nu^{-4}\Vert \Lambda_{t}^{b}\partial_{x}u^{1} \Vert_{L^{2}} \Vert \Lambda_{t}^{b}\partial_{x}b^{2} \Vert_{L^{2}} \Vert \Lambda_{t}^{b}j \Vert_{L^{2}} \\
  & \leq C\nu^{-4}\Vert \Lambda_{t}^{b}\partial_{x}\partial_{y}(-\Delta)^{-1}w_{\neq} \Vert_{L^{2}} \Vert \Lambda_{t}^{b} \partial_{x}\partial_{x}(-\Delta)^{-1}j_{\neq} \Vert_{L^{2}}\Vert \Lambda_{t}^{b}j \Vert_{L^{2}}  \\
  & \leq C\nu^{-4}\Vert \Lambda_{t}^{b} w_{\neq} \Vert_{L^{2}} \Vert \Lambda_{t}^{b} j_{\neq} \Vert_{L^{2}}\Vert \Lambda_{t}^{b}j \Vert_{L^{2}}  \\
  & \leq C\nu^{-4} \Vert (-\Delta)^{-\frac{1}{2}} \Lambda_{t}^{b}j_{\neq} \Vert_{L^{2}}^{\frac{1}{2}} \Vert \nabla \Lambda_{t}^{b}j \Vert_{L^{2}}^{\frac{1}{2}} \Vert (-\Delta)^{-\frac{1}{2}} \Lambda_{t}^{b}w_{\neq} \Vert_{L^{2}}^{\frac{1}{2}} \Vert \nabla \Lambda_{t}^{b}w \Vert_{L^{2}}^{\frac{1}{2}} \Vert \Lambda_{t}^{b}j \Vert_{L^{2}}.
\end{align*}
For $I_{9_{4}} = \langle \Lambda_{t}^{b}(\partial_{x}u^{1}j), \mathcal{M}\Lambda_{t}^{b}j \rangle_{L^{2}}$,
\begin{align*}
  I_{9_{4}} & =  \langle \Lambda_{t}^{b}(j_{0}\partial_{x}u^{1}), \mathcal{M}\Lambda_{t}^{b}j \rangle_{L^{2}} + \langle \Lambda_{t}^{b}(j_{\neq}\partial_{x}u^{1}), \mathcal{M}\Lambda_{t}^{b}j \rangle_{L^{2}}\\
  & =  \langle \Lambda_{t}^{b}(j_{0}\partial_{x}u^{1}), \mathcal{M}\Lambda_{t}^{b}j_{\neq} \rangle_{L^{2}} + \langle \Lambda_{t}^{b}(j_{\neq}\partial_{x}u^{1}), \mathcal{M}\Lambda_{t}^{b}j \rangle_{L^{2}}\\
  & \leq C \nu^{-4}\Vert \Lambda_{t}^{b} j_{0} \Vert_{L^{2}}\Vert (-\Delta)^{-\frac{1}{2}} \Lambda_{t}^{b}j_{\neq} \Vert_{L^{2}}^{\frac{1}{2}} \Vert \nabla \Lambda_{t}^{b}j \Vert_{L^{2}}^{\frac{1}{2}} \Vert (-\Delta)^{-\frac{1}{2}} \Lambda_{t}^{b}w_{\neq} \Vert_{L^{2}}^{\frac{1}{2}} \Vert \nabla \Lambda_{t}^{b}w \Vert_{L^{2}}^{\frac{1}{2}} \\
  & \qquad +  C\nu^{-4} \Vert (-\Delta)^{-\frac{1}{2}} \Lambda_{t}^{b}j_{\neq} \Vert_{L^{2}}^{\frac{1}{2}} \Vert \nabla \Lambda_{t}^{b}j \Vert_{L^{2}}^{\frac{1}{2}} \Vert (-\Delta)^{-\frac{1}{2}} \Lambda_{t}^{b}w_{\neq} \Vert_{L^{2}}^{\frac{1}{2}} \Vert \nabla \Lambda_{t}^{b}w \Vert_{L^{2}}^{\frac{1}{2}} \Vert \Lambda_{t}^{b}j \Vert_{L^{2}}.
\end{align*}
Hence, we have
\begin{equation}\label{9}
\begin{split}
  \vert I_{9} \vert & \leq C\nu^{-4} \Vert (-\Delta)^{-\frac{1}{2}} \Lambda_{t}^{b}j_{\neq} \Vert_{L^{2}}^{\frac{1}{2}} \Vert \nabla \Lambda_{t}^{b}j \Vert_{L^{2}}^{\frac{1}{2}} \Vert (-\Delta)^{-\frac{1}{2}} \Lambda_{t}^{b}w_{\neq} \Vert_{L^{2}}^{\frac{1}{2}} \Vert \nabla \Lambda_{t}^{b}w \Vert_{L^{2}}^{\frac{1}{2}} \Vert \Lambda_{t}^{b}j \Vert_{L^{2}} \\
  & \qquad + C \nu^{-4}\Vert \Lambda_{t}^{b} w_{0} \Vert_{L^{2}}\Vert (-\Delta)^{-\frac{1}{2}} \Lambda_{t}^{b}j_{\neq} \Vert_{L^{2}} \Vert \nabla \Lambda_{t}^{b}j \Vert_{L^{2}}.
\end{split}
\end{equation}
We decompose $I_{10}$ as $I_{10_{1}} + I_{10_{2}}$ with
\[
  I_{10_{1}} = \langle \Lambda_{t}^{b}(\mathbf{u}_{\neq}\cdot \nabla \theta), \vert D_{x} \vert^{\frac{2}{3}}\mathcal{M}\Lambda_{t}^{b}\theta \rangle_{L^{2}},\quad
   I_{10_{2}} = \langle \Lambda_{t}^{b}(\mathbf{u}_{0}\cdot \nabla \theta), \vert D_{x} \vert^{\frac{2}{3}}\mathcal{M}\Lambda_{t}^{b}\theta \rangle_{L^{2}}.
\]
For $I_{10_{1}}$, we have
\begin{align*}
\vert I_{10_{1}} \vert & \leq C \nu^{-4}\Vert \vert D_{x} \vert^{\frac{1}{3}}\Lambda_{t}^{b}(\mathbf{u}_{\neq}\cdot \nabla \theta) \Vert_{L^{2}} \Vert \vert D_{x} \vert^{\frac{1}{3}}\Lambda_{t}^{b}\theta \Vert_{L^{2}} \\ 
& \leq C\nu^{-4} \Vert \Lambda_{t}^{b}\mathbf{u}_{\neq} \Vert_{L^{2}} \Vert \vert D_{x} \vert^{\frac{1}{3}}\Lambda_{t}^{b}\nabla\theta \Vert_{L^{2}}\Vert \vert D_{x} \vert^{\frac{1}{3}}\Lambda_{t}^{b}\theta \Vert_{L^{2}}  \\
& \qquad + C\nu^{-4}\Vert \vert D_{x} \vert^{\frac{1}{3}} \Lambda_{t}^{b}\mathbf{u}_{\neq} \Vert_{L^{2}} \Vert \Lambda_{t}^{b}\nabla\theta \Vert_{L^{2}}\Vert \vert D_{x} \vert^{\frac{1}{3}}\Lambda_{t}^{b}\theta \Vert_{L^{2}} \\
& \leq C\nu^{-4} \Vert (-\Delta)^{-\frac{1}{2}}\Lambda_{t}^{b}w_{\neq} \Vert_{L^{2}} \Vert \vert D_{x} \vert^{\frac{1}{3}}\Lambda_{t}^{b}\nabla\theta \Vert_{L^{2}}\Vert \vert D_{x} \vert^{\frac{1}{3}}\Lambda_{t}^{b}\theta \Vert_{L^{2}} \\
 & \qquad + C\nu^{-4} \Vert \vert D_{x} \vert^{\frac{1}{3}} \Lambda_{t}^{b}w \Vert_{L^{2}} \Vert \Lambda_{t}^{b}\nabla\theta \Vert_{L^{2}}\Vert \vert D_{x} \vert^{\frac{1}{3}}\Lambda_{t}^{b}\theta \Vert_{L^{2}}.
\end{align*}
The estimates for  $I_{10_{2}}$ are the same as those for $I_{1_{2}}$,
\begin{equation*}
\begin{split}
  \vert I_{10_{2}} \vert & \leq C \nu^{-\frac{5}{3}} \Vert \Lambda_{t}^{b}w_{0} \Vert_{L^{2}} \Vert \vert D_{x} \vert^{\frac{2}{3}} \Lambda_{t}^{b}\theta \Vert^{2}_{L^{2}}
   + C\nu^{-5} \Vert \Lambda_{t}^{b}w_{0} \Vert_{L^{2}} \Vert \vert D_{x} \vert^{\frac{1}{3}} \Lambda_{t}^{b}\theta_{\neq} \Vert^{2}_{L^{2}}.
\end{split}
\end{equation*}
Therefore, we deduce that
\begin{align}
  \vert I_{10} \vert & \leq C\nu^{-4}\Vert (-\Delta)^{-\frac{1}{2}}\Lambda_{t}^{b}w_{\neq} \Vert_{L^{2}} \Vert \vert D_{x} \vert^{\frac{1}{3}}\Lambda_{t}^{b}\nabla\theta \Vert_{L^{2}}\Vert \vert D_{x} \vert^{\frac{1}{3}}\Lambda_{t}^{b}\theta \Vert_{L^{2}} \notag\\
   & \qquad + C\nu^{-5} \Vert \Lambda_{t}^{b}w_{0} \Vert_{L^{2}} \Vert (-\Delta)^{-\frac{1}{2}} \vert D_{x} \vert^{\frac{1}{3}}\Lambda_{t}^{b}\theta_{\neq} \Vert_{L^{2}} \Vert \vert D_{x} \vert^{\frac{1}{3}}\nabla \Lambda_{t}^{b}\theta \Vert_{L^{2}} \label{10}\\
     & \qquad + C\nu^{-4}\Vert \vert D_{x} \vert^{\frac{1}{3}} \Lambda_{t}^{b}w \Vert_{L^{2}} \Vert \Lambda_{t}^{b}\nabla\theta \Vert_{L^{2}}\Vert \vert D_{x} \vert^{\frac{1}{3}}\Lambda_{t}^{b}\theta \Vert_{L^{2}}  + C \nu^{-\frac{5}{3}} \Vert \Lambda_{t}^{b}w_{0} \Vert_{L^{2}} \Vert \vert D_{x} \vert^{\frac{2}{3}} \Lambda_{t}^{b}\theta \Vert^{2}_{L^{2}}\notag .
\end{align}
Inserting the upper bounds (\ref{1}){--}(\ref{10}) into the estimates (\ref{bian-1}), (\ref{bian-2}), (\ref{bian-3}) and (\ref{bian-4}), and integrating in time, we obtain
\begin{equation}\label{priori_inequation_theta}
\begin{split}
  & \Vert \Lambda_{t}^{b}\theta \Vert_{L_{t}^{\infty}(L^{2})}^{2} + \nu \Vert \nabla \Lambda_{t}^{b}\theta \Vert_{L_{t}^{2}(L^{2})}^{2} + \frac{1}{4}\nu^{\frac{1}{3}}\Vert \vert D_{x} \vert^{\frac{1}{3}} \Lambda_{t}^{b}\theta \Vert_{L_{t}^{2}(L^{2})}^{2} + \Vert (-\Delta)^{-\frac{1}{2}} \Lambda_{t}^{b}\theta_{\neq} \Vert_{L_{t}^{2}(L^{2})}^{2} \\
  & \leq 2 \Vert \Lambda_{0}^{b}\theta(0) \Vert_{L^{2}}^{2} + C_{1} \nu^{-4} \Vert (-\Delta)^{-\frac{1}{2}}\Lambda_{t}^{b}w_{\neq} \Vert_{L_{t}^{2}(L^{2})}\Vert \nabla \Lambda_{t}^{b}  \theta \Vert_{L^{2}_{t}(L^{2})}\Vert \Lambda_{t}^{b}\theta \Vert_{L^{\infty}_{t}(L^{2})} \\
  & \qquad + C_{1} \nu^{-\frac{5}{3}} \Vert \Lambda_{t}^{b}w \Vert_{L_{t}^{\infty}(L^{2})} \Vert \vert D_{x} \vert^{\frac{1}{3}} \Lambda_{t}^{b}\theta \Vert^{2}_{L_{t}^{2}(L^{2})} \\
  & \qquad + C_{1}\nu^{-5} \Vert \Lambda_{t}^{b}w \Vert_{L_{t}^{\infty}(L^{2})} \Vert (-\Delta)^{-\frac{1}{2}} \Lambda_{t}^{b}\theta_{\neq} \Vert_{L_{t}^{2}(L^{2})} \Vert \nabla \Lambda_{t}^{b}\theta \Vert_{L_{t}^{2}(L^{2})},
\end{split}
\end{equation}
\begin{align}
  & \Vert \Lambda_{t}^{b}w \Vert_{L_{t}^{\infty}(L^{2})}^{2} + \nu \Vert \nabla \Lambda_{t}^{b}w \Vert_{L_{t}^{2}(L^{2})}^{2} + \frac{1}{8}\nu^{\frac{1}{3}}\Vert \vert D_{x} \vert^{\frac{1}{3}} \Lambda_{t}^{b}w \Vert_{L_{t}^{2}(L^{2})}^{2}+ \Vert (-\Delta)^{-\frac{1}{2}} \Lambda_{t}^{b}w_{\neq} \Vert_{L_{t}^{2}(L^{2})}^{2} \notag \\
  & \qquad + \Vert \Lambda_{t}^{b}j \Vert_{L_{t}^{\infty}(L^{2})}^{2}+ \nu \Vert \nabla \Lambda_{t}^{b}j \Vert_{L_{t}^{2}(L^{2})}^{2} + \frac{1}{8}\nu^{\frac{1}{3}}\Vert \vert D_{x} \vert^{\frac{1}{3}} \Lambda_{t}^{b}j \Vert_{L_{t}^{2}(L^{2})}^{2} + \Vert (-\Delta)^{-\frac{1}{2}} \Lambda_{t}^{b}j_{\neq} \Vert_{L_{t}^{2}(L^{2})}^{2} \notag\\
  & \leq 2 \Vert \Lambda_{0}^{b}w(0) \Vert_{L^{2}}^{2} + 2 \Vert \Lambda_{0}^{b}j(0) \Vert_{L^{2}}^{2} + C_{2}\nu^{-\frac{25}{3}}\Vert \vert D_{x} \vert^{\frac{2}{3}} \Lambda_{t}^{b}\theta \Vert_{L_{t}^{2}(L^{2})}^{2}\notag \\
   & \qquad + C_{2} \nu^{-\frac{5}{3}} \Vert \Lambda_{t}^{b}w \Vert_{L_{t}^{\infty}(L^{2})} \Vert \vert D_{x} \vert^{\frac{1}{3}} \Lambda_{t}^{b}w \Vert^{2}_{L_{t}^{2}(L^{2})} 
  + C_{2} \nu^{-\frac{5}{3}} \Vert \Lambda_{t}^{b}w \Vert_{L_{t}^{\infty}(L^{2})} \Vert \vert D_{x} \vert^{\frac{1}{3}} \Lambda_{t}^{b}j \Vert^{2}_{L_{t}^{2}(L^{2})} \notag\\
  & \qquad +   C_{2}\nu^{-4} \Vert (-\Delta)^{-\frac{1}{2}}\Lambda_{t}^{b}w_{\neq} \Vert_{L_{t}^{2}(L^{2})}\Vert \nabla \Lambda_{t}^{b}  w \Vert_{L_{t}^{2}(L^{2})}\Vert \Lambda_{t}^{b}w \Vert_{L_{t}^{\infty}(L^{2})} \notag\\
   & \qquad +  C_{2} \nu^{-5} \Vert \Lambda_{t}^{b}w \Vert_{L_{t}^{\infty}(L^{2})} \Vert (-\Delta)^{-\frac{1}{2}} \Lambda_{t}^{b}w_{\neq} \Vert_{L_{t}^{2}(L^{2})} \Vert \nabla \Lambda_{t}^{b}w \Vert_{L_{t}^{2}(L^{2})}\notag\\
   & \qquad +2 C_{2} \nu^{-4} \Vert (-\Delta)^{-\frac{1}{2}}\Lambda_{t}^{b}j_{\neq} \Vert_{L_{t}^{2}(L^{2})}\Vert \nabla \Lambda_{t}^{b}  j \Vert_{L_{t}^{2}(L^{2})}\Vert \Lambda_{t}^{b}w \Vert_{L_{t}^{\infty}(L^{2})} \notag\\
   & \qquad +2 C_{2} \nu^{-\frac{5}{3}} \Vert \Lambda_{t}^{b}j \Vert_{L_{t}^{\infty}(L^{2})} \Vert \vert D_{x} \vert^{\frac{1}{3}} \Lambda_{t}^{b}j \Vert_{L_{t}^{2}(L^{2})} \Vert \vert D_{x} \vert^{\frac{1}{3}} \Lambda_{t}^{b}w \Vert_{L_{t}^{2}(L^{2})}\notag \\
   & \qquad + 2C_{2}\nu^{-5} \Vert \Lambda_{t}^{b}j \Vert_{L_{t}^{\infty}(L^{2})} \Vert (-\Delta)^{-\frac{1}{2}} \Lambda_{t}^{b}j_{\neq} \Vert_{L_{t}^{2}(L^{2})}^{\frac{1}{2}} \Vert \nabla \Lambda_{t}^{b}j \Vert_{L_{t}^{2}(L^{2})}^{\frac{1}{2}} \label{priori_inequation_w_j} \\
   & \qquad \qquad \times\Vert (-\Delta)^{-\frac{1}{2}} \Lambda_{t}^{b}w_{\neq} \Vert_{L_{t}^{2}(L^{2})}^{\frac{1}{2}} \Vert \nabla \Lambda_{t}^{b}w \Vert_{L_{t}^{2}(L^{2})}^{\frac{1}{2}}\notag \\
   & \qquad +  C_{2} \nu^{-4}\Vert (-\Delta)^{-\frac{1}{2}}\Lambda_{t}^{b}w_{\neq} \Vert_{L_{t}^{2}(L^{2})}\Vert \nabla \Lambda_{t}^{b} j \Vert_{L_{t}^{2}(L^{2})}\Vert \Lambda_{t}^{b}j \Vert_{L_{t}^{\infty}(L^{2})} \notag\\
   & \qquad + C_{2} \nu^{-5} \Vert \Lambda_{t}^{b}w \Vert_{L_{t}^{\infty}(L^{2})} \Vert (-\Delta)^{-\frac{1}{2}} \Lambda_{t}^{b}j_{\neq} \Vert_{L_{t}^{2}(L^{2})} \Vert \nabla \Lambda_{t}^{b}j \Vert_{L_{t}^{2}(L^{2})} \notag 
   \\
   & \qquad + C_{2} \nu^{-4}\Vert (-\Delta)^{-\frac{1}{2}}\Lambda_{t}^{b}j_{\neq} \Vert_{L_{t}^{2}(L^{2})}\Vert \nabla \Lambda_{t}^{b} w \Vert_{L_{t}^{2}(L^{2})}\Vert \Lambda_{t}^{b}j \Vert_{L_{t}^{\infty}(L^{2})} \notag\\
   & \qquad + C_{2}\nu^{-4} \Vert (-\Delta)^{-\frac{1}{2}} \Lambda_{t}^{b}j_{\neq} \Vert_{L^{2}_{t}(L^{2})}^{\frac{1}{2}} \Vert \nabla \Lambda_{t}^{b}j \Vert_{L^{2}_{t}(L^{2})}^{\frac{1}{2}}  \notag\\
   & \qquad \qquad\times \Vert (-\Delta)^{-\frac{1}{2}} \Lambda_{t}^{b}w_{\neq} \Vert_{L^{2}_{t}(L^{2})}^{\frac{1}{2}} \Vert \nabla \Lambda_{t}^{b}w \Vert_{L^{2}_{t}(L^{2})}^{\frac{1}{2}} \Vert \Lambda_{t}^{b}j \Vert_{L^{\infty}_{t}(L^{2})}, \notag 
\end{align}
and
\begin{equation}\label{priori_inequation_dx13_theta}
\begin{split}
  & \Vert \vert D_{x} \vert^{\frac{1}{3}} \Lambda_{t}^{b}\theta \Vert_{L_{t}^{\infty}(L^{2})}^{2} + \nu \Vert \nabla \vert D_{x} \vert^{\frac{1}{3}} \Lambda_{t}^{b}\theta \Vert_{L_{t}^{2}(L^{2})}^{2} + \frac{1}{4}\nu^{\frac{1}{3}}\Vert \vert D_{x} \vert^{\frac{2}{3}} \Lambda_{t}^{b}\theta \Vert_{L_{t}^{2}(L^{2})}^{2} \\
  & \qquad + \Vert (-\Delta)^{-\frac{1}{2}} \vert D_{x} \vert^{\frac{1}{3}} \Lambda_{t}^{b}\theta_{\neq} \Vert_{L_{t}^{2}(L^{2})}^{2} \\
  & \leq 2 \Vert \vert D_{x} \vert^{\frac{1}{3}} \Lambda_{0}^{b}\theta(0) \Vert_{L^{2}}^{2} + C_{3} \nu^{-\frac{5}{3}} \Vert \Lambda_{t}^{b}w \Vert_{L_{t}^{\infty}(L^{2})} \Vert \vert D_{x} \vert^{\frac{2}{3}} \Lambda_{t}^{b}\theta \Vert^{2}_{L_{t}^{2}(L^{2})} \\
  & \qquad + C_{3} \nu^{-4}\Vert (-\Delta)^{-\frac{1}{2}}\Lambda_{t}^{b}w_{\neq} \Vert_{L_{t}^{2}(L^{2})} \Vert \vert D_{x} \vert^{\frac{1}{3}}\Lambda_{t}^{b}\nabla\theta \Vert_{L_{t}^{2}(L^{2})}\Vert \vert D_{x} \vert^{\frac{1}{3}}\Lambda_{t}^{b}\theta \Vert_{L^{\infty}_{t}(L^{2})} \\
  & \qquad + C_{3}\nu^{-4} \Vert \vert D_{x} \vert^{\frac{1}{3}} \Lambda_{t}^{b}w \Vert_{L_{t}^{2}(L^{2})} \Vert \Lambda_{t}^{b}\nabla\theta \Vert_{L_{t}^{2}(L^{2})}\Vert \vert D_{x} \vert^{\frac{1}{3}}\Lambda_{t}^{b}\theta \Vert_{L_{t}^{\infty}(L^{2})} \\
   & \qquad + C_{3} \nu^{-5}\Vert \Lambda_{t}^{b}w \Vert_{L_{t}^{\infty}(L^{2})} \Vert (-\Delta)^{-\frac{1}{2}} \vert D_{x} \vert^{\frac{1}{3}}\Lambda_{t}^{b}\theta_{\neq} \Vert_{L_{t}^{2}(L^{2})} \Vert \vert D_{x} \vert^{\frac{1}{3}}\nabla \Lambda_{t}^{b}\theta \Vert_{L^{2}_{t}(L^{2})}.
\end{split}
\end{equation}

The priori bounds in (\ref{priori_inequation_theta}), (\ref{priori_inequation_w_j}) and (\ref{priori_inequation_dx13_theta}) allow us to prove Theorem \ref{nolinear_theorem} through the bootstrap argument. We recall the assumptions on the initial data $(w(0),j(0),\theta(0))$,
\begin{equation*}
  \Vert \theta(0) \Vert_{H^{b}} \leq \varepsilon \nu^{\alpha},\qquad \Vert (w(0),j(0)) \Vert_{H^{b}} \leq \varepsilon \nu^{\beta},\qquad \Vert \vert D_{x} \vert^{\frac{1}{3}} \theta(0) \Vert_{H^{b}} \leq \varepsilon \nu^{\delta},
\end{equation*}
where $\varepsilon > 0$ is sufficiently small and
\begin{equation}\label{initial_condition_}
  \beta \geq \frac{11}{2}, \qquad \delta \geq \beta + \frac{13}{3}, \qquad \alpha \geq \delta - \beta + \frac{14}{3}.
\end{equation}
To apply the bootstrap argument, we make the ansatz that, for $T\leq \infty$, the solution of (\ref{nonlinear_stability_system_lambdatb}) obeys
\begin{align}
 \begin{split}\label{ansatz1_theta}
  & \Vert \Lambda_{t}^{b}\theta \Vert_{L_{t}^{\infty}([0,T])(L^{2})} + \nu^{\frac{1}{2}}\Vert \nabla \Lambda_{t}^{b}\theta \Vert_{L_{t}^{2}([0,T])(L^{2})}  +  \nu^{\frac{1}{6}}\Vert \vert D_{x} \vert^{\frac{1}{3}} \Lambda_{t}^{b}\theta \Vert_{L_{t}^{2}([0,T])(L^{2})}  \\ 
  &\qquad + \Vert (-\Delta)^{-\frac{1}{2}}\Lambda_{t}^{b}\theta_{\neq} \Vert_{L_{t}^{2}([0,T])(L^{2})}  \leq C\varepsilon  \nu^{\alpha},
\end{split} 
\end{align}
\begin{align}
\begin{split}\label{ansatz2_w_and_j}
  &\Vert \Lambda_{t}^{b}w \Vert_{L_{t}^{\infty}([0,T])(L^{2})} + \nu^{\frac{1}{2}}\Vert \nabla \Lambda_{t}^{b}w \Vert_{L_{t}^{2}([0,T])(L^{2})}  + \nu^{\frac{1}{6}}\Vert \vert D_{x} \vert^{\frac{1}{3}} \Lambda_{t}^{b}w \Vert_{L_{t}^{2}([0,T])(L^{2})}    \\
  & \qquad  + \Vert \Lambda_{t}^{b}j \Vert_{L_{t}^{\infty}([0,T])(L^{2})} + \nu^{\frac{1}{2}}\Vert \nabla \Lambda_{t}^{b}j \Vert_{L_{t}^{2}([0,T])(L^{2})}+\nu^{\frac{1}{6}}\Vert \vert D_{x} \vert^{\frac{1}{3}} \Lambda_{t}^{b}j \Vert_{L_{t}^{2}([0,T])(L^{2})}  \\
  &\qquad+ \Vert (-\Delta)^{-\frac{1}{2}}\Lambda_{t}^{b}w_{\neq} \Vert_{L_{t}^{2}([0,T])(L^{2})}  + \Vert (-\Delta)^{-\frac{1}{2}}\Lambda_{t}^{b}j_{\neq} \Vert_{L_{t}^{2}([0,T])(L^{2})} \leq C\varepsilon   \nu^{\beta},
\end{split} 
\end{align}
\begin{align}
\begin{split}\label{ansatz4_dx13_theta}
  &\Vert \vert D_{x} \vert^{\frac{1}{3}} \Lambda_{t}^{b}\theta \Vert_{L_{t}^{\infty}([0,T])(L^{2})} + \nu^{\frac{1}{2}}\Vert \nabla \vert D_{x} \vert^{\frac{1}{3}} \Lambda_{t}^{b}\theta \Vert_{L_{t}^{2}([0,T])(L^{2})} + \nu^{\frac{1}{6}}\Vert \vert D_{x} \vert^{\frac{2}{3}} \Lambda_{t}^{b}\theta \Vert_{L_{t}^{2}([0,T])(L^{2})}  \\
 &\qquad + \Vert (-\Delta)^{-\frac{1}{2}}\vert D_{x} \vert^{\frac{1}{3}}\Lambda_{t}^{b}\theta_{\neq} \Vert_{L_{t}^{2}([0,T])(L^{2})}  \leq \tilde{C} \varepsilon   \nu^{\delta}.
\end{split}
\end{align}
We then show that (\ref{ansatz1_theta}), (\ref{ansatz2_w_and_j}), and (\ref{ansatz4_dx13_theta}) actually hold with $C$ replaced by $C/2$ and $\tilde{C}$ by $\tilde{C}/2$. In fact, if we insert the initial condition and the ansatz in priori, we find
\begin{align*}
  & \Vert \Lambda_{t}^{b}\theta \Vert_{L_{t}^{\infty}([0,T])(L^{2})}^{2} + \nu \Vert \nabla \Lambda_{t}^{b}\theta \Vert_{L_{t}^{2}([0,T])(L^{2})}^{2} + \frac{1}{4}\nu^{\frac{1}{3}}\Vert \vert D_{x} \vert^{\frac{1}{3}} \Lambda_{t}^{b}\theta \Vert_{L_{t}^{2}([0,T])(L^{2})}^{2} \\
  & \quad + \Vert (-\Delta)^{-\frac{1}{2}} \Lambda_{t}^{b}\theta_{\neq} \Vert_{L_{t}^{2}([0,T])(L^{2})}^{2} \\
  & \leq 2\varepsilon^{2}\nu^{2\alpha} + C_{1}C^{3}\varepsilon^{3}(\nu^{2\alpha + \beta - \frac{9}{2}} + \nu^{2\alpha + \beta-2} + \nu^{2\alpha + \beta - \frac{11}{2}}),
  \end{align*}
  \begin{align*}
  & \Vert \Lambda_{t}^{b}w \Vert_{L_{t}^{\infty}([0,T])(L^{2})}^{2} + \nu \Vert \nabla \Lambda_{t}^{b}w \Vert_{L_{t}^{2}([0,T])(L^{2})}^{2} + \frac{1}{8}\nu^{\frac{1}{3}}\Vert \vert D_{x} \vert^{\frac{1}{3}} \Lambda_{t}^{b}w \Vert_{L_{t}^{2}([0,T])(L^{2})}^{2} \\
  & \quad+ \Vert \Lambda_{t}^{b}j \Vert_{L_{t}^{\infty}([0,T])(L^{2})}^{2} + \nu \Vert \nabla \Lambda_{t}^{b}j \Vert_{L_{t}^{2}([0,T])(L^{2})}^{2} + \frac{1}{8}\nu^{\frac{1}{3}}\Vert \vert D_{x} \vert^{\frac{1}{3}} \Lambda_{t}^{b}j \Vert_{L_{t}^{2}([0,T])(L^{2})}^{2} \\
  &  \quad + \Vert (-\Delta)^{-\frac{1}{2}} \Lambda_{t}^{b}w_{\neq} \Vert_{L_{t}^{2}([0,T])(L^{2})}^{2}  + \Vert (-\Delta)^{-\frac{1}{2}} \Lambda_{t}^{b}j_{\neq} \Vert_{L_{t}^{2}([0,T])(L^{2})}^{2} \\
  &  \leq 2\varepsilon^{2}\nu^{2\beta} + 2\varepsilon^{2}\nu^{2\beta} + C_{2} \tilde{C}^{2} \varepsilon^{2} \nu^{2\delta - \frac{26}{3}} + C_{2}C^{3}\varepsilon^{3}\nu^{3\beta - \frac{9}{2}} + C_{2}C^{3}\varepsilon^{3}\nu^{3\beta-2} + C_{2}C^{3}\varepsilon^{3}\nu^{3\beta - \frac{11}{2}} \\
  & \quad + C_{2}C^{3}\varepsilon^{3}\nu^{3\beta - \frac{9}{2}} 
  + 2C_{2}C^{3}\varepsilon^{3}\nu^{3\beta-2} 
  + 2C_{2}C^{3}\varepsilon^{3}\nu^{3\beta - \frac{11}{2}}
   +2 C_{2}C^{3}\varepsilon^{3}\nu^{3\beta - \frac{9}{2}} 
   + C_{2}C^{3}\varepsilon^{3}\nu^{3\beta-2} \\
  &\quad  + C_{2}C^{3}\varepsilon^{3}\nu^{3\beta - \frac{11}{2}} 
  +  C_{2}C^{3}\varepsilon^{3}\nu^{3\beta - \frac{9}{2}} 
   + C_{2}C^{3}\varepsilon^{3}\nu^{3\beta- \frac{9}{2}} \\
 &  \leq 4\varepsilon^{2}\nu^{2\beta} + C_{2} \tilde{C}^{2} \varepsilon^{2} \nu^{2\delta - \frac{26}{3}} + C_{2}C^{3}\varepsilon^{3}(6\nu^{3\beta - \frac{9}{2}} + 4\nu^{3\beta - 2}+ 4\nu^{3\beta - \frac{11}{2}}),
\end{align*}
\begin{align*}
  & \Vert \vert D_{x} \vert^{\frac{1}{3}} \Lambda_{t}^{b}\theta \Vert_{L_{t}^{\infty}([0,T])(L^{2})}^{2} + \nu \Vert \nabla \vert D_{x} \vert^{\frac{1}{3}} \Lambda_{t}^{b}\theta \Vert_{L_{t}^{2}([0,T])(L^{2})}^{2} + \frac{1}{4}\nu^{\frac{1}{3}}\Vert \vert D_{x} \vert^{\frac{2}{3}} \Lambda_{t}^{b}\theta \Vert_{L_{t}^{2}([0,T])(L^{2})}^{2} \\
  & \quad+ \Vert (-\Delta)^{-\frac{1}{2}} \vert D_{x} \vert^{\frac{1}{3}} \Lambda_{t}^{b}\theta_{\neq} \Vert_{L_{t}^{2}([0,T])(L^{2})}^{2} \\
  & \leq 2\varepsilon^{2}\nu^{2\delta} + C_{3}C \tilde{C}^{2} \varepsilon^{3} \nu^{2\delta + \beta - \frac{9}{2}} + C_{3}C^{2} \tilde{C} \varepsilon^{3} \nu^{\alpha + \delta + \beta - \frac{14}{3}} \\
  & \quad + C_{3}C \tilde{C}^{2}\varepsilon^{3} \nu^{2\delta + \beta - 2} + C_{3}C \tilde{C}^{2} \varepsilon^{3}\nu^{2\delta + \beta - \frac{11}{2}} \\
  & \leq 2\varepsilon^{2}\nu^{2\delta} + C_{3}C \tilde{C} \varepsilon^{3} (\tilde{C}\nu^{2\delta + \beta - \frac{9}{2}} +  C \nu^{\alpha + \delta + \beta - \frac{14}{3}} + \tilde{C} \nu^{2\delta + \beta-2} + \tilde{C} \nu^{2\delta + \beta - \frac{11}{2}}),
\end{align*}
which implies that
\begin{align*}
  & \Vert \Lambda_{t}^{b}\theta \Vert_{L_{t}^{\infty}([0,T])(L^{2})} + \nu^{\frac{1}{2}} \Vert \nabla \Lambda_{t}^{b}\theta \Vert_{L_{t}^{2}([0,T])(L^{2})} + \nu^{\frac{1}{6}}\Vert \vert D_{x} \vert^{\frac{1}{3}} \Lambda_{t}^{b}\theta \Vert_{L_{t}^{2}([0,T])(L^{2})} \\
  & \quad + \Vert (-\Delta)^{-\frac{1}{2}} \Lambda_{t}^{b}\theta_{\neq} \Vert_{L_{t}^{2}([0,T])(L^{2})} \\
  &  \leq 4 \varepsilon \nu^{\alpha} + 3 C_{1}^{\frac{1}{2}} C^{\frac{3}{2}} \varepsilon^{\frac{3}{2}} (\nu^{\alpha + \frac{1}{2}\beta - \frac{9}{4}} + \nu^{\alpha + \frac{1}{2}\beta-1} + \nu^{\alpha + \frac{1}{2}\beta -\frac{11}{4}}),
\end{align*}
\begin{align*}
  & \Vert \Lambda_{t}^{b}w \Vert_{L_{t}^{\infty}([0,T])(L^{2})} + \nu^{\frac{1}{2}} \Vert \nabla \Lambda_{t}^{b}w \Vert_{L_{t}^{2}([0,T])(L^{2})} + \nu^{\frac{1}{6}}\Vert \vert D_{x} \vert^{\frac{1}{3}} \Lambda_{t}^{b}w \Vert_{L_{t}^{2}([0,T])(L^{2})} \\
  & \quad+ \Vert \Lambda_{t}^{b}j \Vert_{L_{t}^{\infty}([0,T])(L^{2})}^{2} + \nu^{\frac{1}{2}} \Vert \nabla \Lambda_{t}^{b}j \Vert_{L_{t}^{2}([0,T])(L^{2})}^{2} + \nu^{\frac{1}{6}}\Vert \vert D_{x} \vert^{\frac{1}{3}} \Lambda_{t}^{b}j \Vert_{L_{t}^{2}([0,T])(L^{2})}^{2} \\
  & \quad   + \Vert (-\Delta)^{-\frac{1}{2}} \Lambda_{t}^{b}w_{\neq} \Vert_{L_{t}^{2}([0,T])(L^{2})} 
   + \Vert (-\Delta)^{-\frac{1}{2}} \Lambda_{t}^{b}j_{\neq} \Vert_{L_{t}^{2}([0,T])(L^{2})}^{2} \\
 &  \leq 10 \varepsilon \nu^{\beta} + 5C_{2}^{\frac{1}{2}} \tilde{C} \varepsilon \nu^{\delta - \frac{13}{3}} + 12 C_{2}^{\frac{1}{2}}C^{\frac{3}{2}}\varepsilon^{\frac{3}{2}}(\nu^{\frac{3}{2}\beta - \frac{9}{4}} + \nu^{\frac{3}{2}\beta - 1} + \nu^{\frac{3}{2}\beta-\frac{11}{4}}),
\end{align*}
\begin{align*}
  & \Vert \vert D_{x} \vert^{\frac{1}{3}} \Lambda_{t}^{b}\theta \Vert_{L_{t}^{\infty}([0,T])(L^{2})} + \nu^{\frac{1}{2}} \Vert \nabla \vert D_{x} \vert^{\frac{1}{3}} \Lambda_{t}^{b}\theta \Vert_{L_{t}^{2}([0,T])(L^{2})} + \nu^{\frac{1}{6}}\Vert \vert D_{x} \vert^{\frac{2}{3}} \Lambda_{t}^{b}\theta \Vert_{L_{t}^{2}([0,T])(L^{2})} \\
  &\quad + \Vert (-\Delta)^{-\frac{1}{2}} \vert D_{x} \vert^{\frac{1}{3}} \Lambda_{t}^{b}\theta_{\neq} \Vert_{L_{t}^{2}([0,T])(L^{2})} \\
  &  \leq 4 \varepsilon \nu^{\delta}  + 3C_{3}^{\frac{1}{2}}C^{\frac{1}{2}}\tilde{C}^{\frac{1}{2}}\varepsilon^{\frac{3}{2}}(\tilde{C}^{\frac{1}{2}}\nu^{\delta+\frac{1}{2}\beta-\frac{9}{4}}+ C^{\frac{1}{2}}\nu^{\frac{1}{2}\alpha + \frac{1}{2}\delta+\frac{1}{2}\beta-\frac{7}{3}} + \tilde{C}^{\frac{1}{2}}\nu^{\delta+\frac{1}{2}\beta-1} + \tilde{C}^{\frac{1}{2}}\nu^{\delta+\frac{1}{2}\beta - \frac{11}{4}}).
\end{align*}
If we invoke (\ref{initial_condition_}) and choose
\[
  \tilde{C} \geq 32, C \geq 80 ,\ C \geq 40 C_{2}^{\frac{1}{2}}\tilde{C},\  \varepsilon \leq \min\lbrace{\frac{1}{288^{2} C C_{2}}, \frac{\tilde{C}}{24^{2}C^{2}C_{3}}, \frac{1}{72^{2} C C_{3}}, \frac{1}{36^{2}C_{1}C} \rbrace},
\]
then the inequalities hold with $C$ replaced by $C/2$ and $\tilde{C}$ replaced by $\tilde{C}/2$. Hence, the ansatz is true for all $T \geq 0$. We redefine $C$ as the $\max\{C,\tilde{C}\}$ and this completes the proof. 

\bigskip
\noindent {\bf Acknowledgment.}
This work is supported by NSFC under the contracts 11871005 and 11771041.

\end{document}